\documentclass[10pt]{article}
\usepackage[utf8]{inputenc}
\usepackage[a4paper,top=2.5cm,bottom=2cm,left=2cm,right=2cm,marginparwidth=1.5cm]{geometry}
\usepackage{amsmath}
\usepackage{amssymb}
\numberwithin{equation}{section}
\usepackage[english]{babel}
\usepackage{amsthm}
\usepackage{bbm}
\usepackage{amsfonts}
\usepackage{comment}
\usepackage{multirow}
\usepackage{mathrsfs}
\usepackage{hyperref}
\usepackage{titlesec}
\usepackage{bbm}
\usepackage[numbers,sort]{natbib}
\hypersetup{
    colorlinks=true,
    linkcolor=blue,
    citecolor=red,
    filecolor=magenta,      
    urlcolor=black,
    pdftitle={MLResonances},
    pdfpagemode=FullScreen
    }
\usepackage{graphicx, tikz}
\usepackage{subcaption}
\usepackage{float}
\usepackage{xcolor}
\usepackage{bbm}
\usepackage[format=plain,
            font=it]{caption}

\newcommand{\upd}{\mathrm{d}}

       \def\w{\omega }

\title{Asymptotics-guided learning and symbolic regression for dispersive resonances}
\author{Konstantinos Alexopoulos\thanks{CMAP, CNRS, Ecole polytechnique, Institut Polytechnique de Paris, 91120 Palaiseau, France} \and Josselin Garnier \thanks{CMAP, CNRS, Ecole polytechnique, Institut Polytechnique de Paris, 91120 Palaiseau, France}}
\date{}

\begin{document}

\maketitle

\begin{abstract}
    We study resonance prediction in dispersive media, formulated as nonlinear spectral problems for volume integral operators. The main idea is to use asymptotic analysis not only as a baseline approximation, but also as a guide for constructing predictive correction models. We learn the residual between asymptotic and reference resonances using features suggested by the subwavelength expansion, including the logarithmic scales specific to two dimensions. The resulting corrections substantially improve single-resonator and dimer predictions, and symbolic regression produces compact formulas for the learned residual. The results show that asymptotic analysis can be used not only to approximate resonances, but also to design the feature space in which data-driven corrections become accurate, low-dimensional, and interpretable.
\end{abstract}

\section{Introduction}

Resonance phenomena in wave propagation play a central role in a wide range of applications, including photonics, metamaterials, and electromagnetic device design. In dispersive media, the computation of resonant frequencies leads to nonlinear spectral problems associated with volume integral formulations, typically derived from the Lippmann--Schwinger equation \cite{ADFMS}. These problems are particularly important in the analysis of subwavelength resonators and highly dispersive systems, where the material parameters depend strongly on frequency.

A prominent example arises in the study of halide perovskite resonators \cite{AM}, whose dispersive permittivity \cite{MFTHBPZK} enables strong light-matter interaction and supports a wide range of applications, including optical sensing \cite{GPLLZZSQKJ}, photovoltaics \cite{S}, and light-emitting devices \cite{WKYZZPLBHL}. More broadly, resonant structures underpin many phenomena in metamaterials, such as negative effective parameters \cite{pendry2000negative, veselago1968electrodynamics}, cloaking \cite{milton2006cloaking, craster2012acoustic}, and bio-inspired optical effects \cite{sun2013structural, zhao2012bio, GPLLZZSQKJ, lee2018bioinspired}. Accurate and efficient prediction of resonances is therefore essential for both mathematical analysis and design.

From a mathematical perspective, asymptotic analysis provides a powerful framework for approximating resonances in subwavelength regimes. Expansions derived from the spectral properties of volume integral operators yield explicit formulas that capture the leading-order behaviour of resonant frequencies \cite{AD, alexopoulos2022mathematical}. These approaches have been successfully applied to both single resonators and coupled systems, where interaction effects lead to mode splitting and branch-dependent resonance behaviour. However, their accuracy deteriorates away from the strict asymptotic regime, and higher-order corrections are often difficult to compute explicitly, especially in the presence of coupling or complex geometries.

In parallel, machine learning techniques have become important tools in the physical sciences \cite{carleo2019machine}, with applications ranging from operator learning and partial differential equations \cite{tsakyridis2024photonic} to inverse design in photonics \cite{ma2021deep, wiecha2021deep, gao2021deep, yun2022deep, hegde2020deep}. In nanophotonics, deep learning has been widely used for the design and optimization of structures with prescribed optical properties. However, many such approaches learn the full input-output map directly from data and rely on high-capacity models, which may require large datasets and can be difficult to interpret.

In this work, we take a different approach. Rather than replacing the asymptotic model, we use it as the leading-order approximation and learn only the residual error between this approximation and accurate reference resonances. The learning task is therefore not to rediscover the full resonance map, but to correct the missing higher-order terms. This residual is modelled using physics-informed features derived from the asymptotic expansion, logarithmic frequency-dependent scales suggested by the two-dimensional subwavelength expansions. This leads to a low-dimensional correction model that remains tied to the underlying mathematical structure.

A second objective of the paper is to obtain interpretable correction formulas. For this purpose, we use symbolic regression, which searches for explicit analytic expressions fitting the data. Symbolic regression has recently received significant attention as a tool for interpretable scientific machine learning, notably through PySR and SymbolicRegression.jl \cite{cranmer2023interpretable}. 
It has also begun to appear in photonics and optical modelling, including applications to topological photonic bands \cite{ghorashi2026symbolic}, optical responses of biological and bio-inspired structures 
\cite{sierra2024symbolic, sierra2025modeling}, and broader discovery and optimization problems in physics and topological photonics \cite{kim2023novel}. Here, symbolic regression is used not as a black-box discovery method, but as an asymptotically guided tool for extracting compact formulas for the learned residual.

The main contributions of this work are as follows. First, we introduce an asymptotics-guided residual learning framework for nonlinear resonance problems in dispersive media, in which machine learning is used to correct, rather than replace, asymptotic or reduced-order models. Second, we construct physics-informed features based on the known structure of the subwavelength expansion, including logarithmic next-order terms, and show through ablation studies that these features are essential for accurate correction. Third, we use symbolic regression to obtain explicit correction formulas for both a single resonator and a coupled dimer, leading to interpretable residual models. Fourth, we show that assisting symbolic regression with a proxy for the next asymptotic correction improves or stabilizes the learned formulas, supporting the interpretation that the residual is governed by the expected asymptotic scales. Finally, we investigate the effect of noisy training resonances and show that the proposed framework remains robust to moderate noise, especially when the number of training samples is increased.

The paper is organized as follows. In Section~\ref{sec:Mathematical setting}, we introduce the mathematical setting and formulate the resonance problem using volume integral operators. Section~\ref{sec:Asymptotic modelling} recalls the asymptotic and reduced modelling framework, including the next-order structures that motivate the feature design. In Section~\ref{sec:Learning framework}, we present the residual learning method based on physics-informed Ridge regression. Section~\ref{sec:symbolic_regression} introduces the symbolic-regression correction and the assisted symbolic-regression strategy. Section~\ref{sec:Numerical experiments} presents the numerical experiments for single-resonator and dimer configurations, including ablation studies, assisted symbolic-regression tests, and noise robustness experiments. Finally, Section~\ref{sec:Conclusion and perspectives} discusses the implications, limitations, and possible extensions of the approach.

\section{Mathematical setting}\label{sec:Mathematical setting}

\subsection{Governing equations, non-dimensionalisation and parametric setting}

We consider time-harmonic wave propagation in a two-dimensional setting, governed by the Helmholtz equation
\begin{equation}
    \nabla \cdot \left( \frac{1}{\mu(x)} \nabla u(x) \right)
    + \omega^2 \varepsilon(x,\omega) u(x) = 0,
    \qquad x \in \mathbb{R}^2,
\end{equation}
where $\omega \in \mathbb{C}$ denotes the angular frequency, $\varepsilon(x,\omega)$ is the possibly dispersive and complex-valued permittivity, and $\mu(x)$ is the permeability. The background medium is homogeneous, with constant parameters $(\varepsilon_0,\mu_0)$, and the resonators occupy a bounded domain $D\subset\mathbb{R}^2$. The material parameters are
\begin{equation}
\varepsilon(x,\omega) =
    \begin{cases}
        \varepsilon_r(\omega), & x \in D, \\
        \varepsilon_0, & x \in \mathbb{R}^2 \setminus D,
    \end{cases}
    \qquad
    \mu(x)\equiv\mu_0 .
\end{equation}
The resonator permittivity is described by a Lorentz-type dispersive law,
\begin{equation}\label{eq:lorentz_model}
    \varepsilon_r(\omega)
    =
    \varepsilon_0
    +
    \frac{\alpha}
    {\beta-\omega^2-i\gamma\omega+\eta k_0^2},
    \qquad
    k_0=\omega\sqrt{\varepsilon_0\mu_0},
\end{equation}
where the material parameters \(\alpha,\beta,\gamma,\eta\) are real, with \(\gamma\geq 0\) describing damping. The frequency \(\omega\in\mathbb{C}\) is complex in the resonance problem and therefore the background wavenumber \(k_0\) may also be complex. This model captures both dispersion and damping effects and is commonly used in the study of perovskite-based resonators, as described in \cite{AD}. Since $k_0^2=\omega^2\varepsilon_0\mu_0$, the denominator in \eqref{eq:lorentz_model} may equivalently be viewed as a quadratic polynomial in $\omega$. We keep the form involving $k_0^2$ because it matches the notation used in the asymptotic expansion of the volume integral operator. 

In order to avoid dimensional ambiguities, we work throughout with non-dimensional variables. We choose the reference frequency and length scales
\begin{equation}
    \omega_*=\sqrt{\beta},
    \qquad
    L_*=\frac{1}{\omega_*\sqrt{\varepsilon_0\mu_0}},
\end{equation}
and define
\begin{equation}
    \widehat{x}=\frac{x}{L_*},
    \qquad
    \widehat{\omega}=\frac{\omega}{\omega_*},
    \qquad
    \widehat{k}_0=L_* k_0,
    \qquad
    \widehat{\gamma}=\frac{\gamma}{\omega_*}.
\end{equation}
The eigenvalues of the volume integral operator are scaled by
\begin{equation}
    \widehat{\lambda}_j=\frac{\lambda_j}{L_*^2}.
\end{equation}
After this non-dimensionalisation, the resonance condition keeps the same form, but all quantities appearing in it are dimensionless. For readability, we drop the hats in the remainder of the paper. Thus $\omega$, $k_0$, $\gamma$, and $\lambda_j$ should be understood below as non-dimensional quantities, and the numerical constants appearing in the symbolic-regression formulas are dimensionless fitted coefficients.

We now specify the parametric family considered in the numerical experiments. The reference single-resonator geometry is fixed and is taken to be the unit disk
$B=\{x\in\mathbb{R}^2:\ |x|<1\}$.
The physical resonator is obtained by the subwavelength rescaling
\begin{equation}
    D_\delta=\delta B,
    \qquad
    0<\delta\ll1.
\end{equation}
For the dimer, the reference geometry depends on a dimensionless separation parameter $\kappa>0$ and is given by
$B_\kappa=B_1^\kappa\cup B_2^\kappa$,
where $B_1^\kappa$ and $B_2^\kappa$ are two unit disks whose centres are located at
$\left(-1-\frac{\kappa}{2},0\right)$ and $
    \left(1+\frac{\kappa}{2},0\right)$.
Thus $\kappa$ represents the gap between the two disks in the reference configuration, and the scaled dimer is
\begin{equation}
    D_{\delta,\kappa}=\delta B_\kappa .
\end{equation}

In the numerical experiments, we denote by \(\omega_{\mathrm{ref}}\) the reference resonance computed from the high-fidelity numerical model introduced below. In the single-resonator case this is one complex resonance, while in the dimer case the two resonance branches are denoted by \(\omega_{1,\mathrm{ref}}\) and \(\omega_{2,\mathrm{ref}}\).

In the learning experiments, the geometry of the reference resonator is fixed. For the single-resonator problem, the only varying parameters are the size parameter $\delta$ and the damping parameter $\gamma$. Thus the learned correction is constructed over the parametric map
\begin{align*}
    (\delta,\gamma) \longmapsto \omega_{\mathrm{ref}}.
\end{align*}
For the dimer, the varying parameters are the size parameter $\delta$, the damping parameter $\gamma$, and the dimensionless separation parameter $\kappa$, giving the map
\begin{align*}
    (\delta,\gamma,\kappa) \longmapsto \big(\omega_{1,\mathrm{ref}},\omega_{2,\mathrm{ref}}\big).
\end{align*}
The material parameters $\alpha,\beta,\eta$ are fixed and assumed to be known. The damping parameter $\gamma$ is varied because it directly controls the imaginary part of the dispersive response and therefore the resonance linewidth. The present work should therefore be understood as a residual correction study for a fixed dispersive material family, not as a learning problem over arbitrary materials.

This choice reflects the physical role of the different parameters in the Lorentz model. The parameters $\alpha$, $\beta$, and $\eta$ determine the underlying dispersive material law: $\alpha$ fixes the oscillator strength, $\beta$ fixes the characteristic material frequency, and $\eta$ controls the strength of the spatial-dispersion correction. In the present study, these quantities are treated as known material constants. By contrast, $\gamma$ controls the damping, and therefore directly affects the imaginary part and linewidth of the resonances. It is also a natural parameter through which losses, fabrication variability, or material quality can change. We therefore vary $\gamma$ in order to test whether the learned correction remains stable across different damping regimes, while keeping the underlying dispersive material family fixed.

Finally, the scattered field is required to satisfy the outgoing Sommerfeld radiation condition. Namely, if $u^{\mathrm{inc}}$ denotes the incident field and $k_0=\omega\sqrt{\varepsilon_0\mu_0}$ is the background wavenumber, then
\begin{equation}\label{Sommerfeld}
    \lim_{|x|\to\infty}
    |x|^{\frac{1}{2}}
    \left(
        \frac{\partial}{\partial |x|}
        -ik_0
    \right)
    \big(u(x)-u^{\mathrm{inc}}(x)\big)
    =0.
\end{equation}

\subsection{Volume integral formulation}

Let us introduce the contrast function
\begin{equation}
    \xi(\omega) = \mu_0 \left( \varepsilon_r(\omega) - \varepsilon_0 \right).
\end{equation}
The total field $u$ can be represented via the Lippmann–Schwinger equation
\begin{equation}
    u(x) = u^{\mathrm{inc}}(x) - \omega^2 \int_D G(x,y;\omega)\, \xi(\omega)\, u(y)\, dy,
\end{equation}
where $G(\cdot,\cdot;\omega)$ is the free-space Green's function of the background medium. In two dimensions, the Green's function admits the representation
\begin{equation}
    G(x,y;\omega) = -\frac{1}{2\pi} \log|x-y| + R(x,y;\omega),
\end{equation}
where $R$ is a smooth remainder term. This logarithmic singularity is responsible for the characteristic logarithmic terms that appear in the subwavelength asymptotics.
Defining the volume integral operator
\begin{equation}
    \mathcal{K}^{\w}_D[u](x) = \int_D G(x,y;\w)\, u(y)\, dy,
    \label{eq:defKD}
\end{equation}
the Lippmann-Schwinger equation can be written in operator form as
$
    u + \omega^2 \xi(\omega) \mathcal{K}^{\w}_D[u] = u^{\mathrm{inc}}
$.
In the subwavelength regime considered below, asymptotic approximations of
$\mathcal{K}_D^\omega$ yield reduced resonance conditions that can be solved efficiently.

\subsection{Resonance problem}

Resonances are defined as complex frequencies $\omega$ for which the corresponding homogeneous problem admits a non-trivial outgoing solution. Equivalently, they are the values of $\omega\in\mathbb{C}$ for which
\begin{equation}
    u + \omega^2 \xi(\omega)\mathcal{K}_D^\omega[u] = 0
\end{equation}
has a non-zero solution. Thus the operator
$
    \mathrm{Id} + \omega^2 \xi(\omega)\mathcal{K}_D^\omega
$
is not invertible. This gives a nonlinear spectral problem, since both the material contrast and the integral operator depend on the frequency.

If $\lambda_j(\omega)$ denotes an eigenvalue of $\mathcal{K}_D^\omega$, then a formal resonance condition is
\begin{equation}
\label{eq:resonance_condition}
    1 + \omega^2 \xi(\omega)\lambda_j(\omega) = 0 .
\end{equation}
In the asymptotic regime, one often replaces $\lambda_j(\omega)$ by an approximation obtained from the leading terms of the volume integral operator. This leads to explicit or low-dimensional approximations of the resonant frequencies. 

The corresponding asymptotic prediction is denoted by \(\omega_{\mathrm{asymp}}\), and the learned residual is \(\Delta \omega = \omega_{\mathrm{ref}}-\omega_{\mathrm{asymp}}\).

\subsection{Coupled resonator systems}

We also consider systems of multiple resonators. In particular, for a pair of identical resonators, or dimer, we use the scaled configuration
$
    D_{\delta,\kappa}=\delta B_\kappa
$,
where $B_\kappa$ is the reference dimer introduced above. Thus $\kappa$ controls the separation in the reference geometry, while the physical configuration is obtained by the common scaling factor $\delta$. For notational simplicity, we write $D$ instead of $D_{\delta,\kappa}$ and $B$ instead of $B_\kappa$ when no confusion is possible.

In this setting, the operator $\mathcal{K}_D^\omega$ has a natural block structure,
\begin{equation}
    \mathcal{K}_D^\omega =
    \begin{pmatrix}
        \mathcal{K}_{D_1D_1}^\omega & \mathcal{K}_{D_1D_2}^\omega \\
        \mathcal{K}_{D_2D_1}^\omega & \mathcal{K}_{D_2D_2}^\omega
    \end{pmatrix},
    \label{eq:defKDomegadimer}
\end{equation}
where the diagonal blocks describe self-interactions and the off-diagonal blocks describe coupling between the two resonators. This coupling splits the leading resonance into two branches, associated with symmetric and antisymmetric modes. The size of the splitting depends on the separation distance $\kappa$, as described in \cite{alexopoulos2022mathematical}.

Consequently, the resonance condition \eqref{eq:resonance_condition} gives two distinct resonant frequencies, which we denote by
\begin{align*}
    \omega_1 \quad \text{and} \quad \omega_2 .
\end{align*}
Their asymptotic approximations will be denoted by
\begin{align*}
    \omega_{1,\mathrm{asymp}} \quad \text{and} \quad \omega_{2,\mathrm{asymp}}
\end{align*}
and the corresponding residuals by
\begin{align*}
    \Delta\omega_j =
    \omega_{j,\mathrm{ref}}-\omega_{j,\mathrm{asymp}},
    \qquad j=1,2.
\end{align*}
The accurate approximation of these residuals is more delicate than in the single-resonator case, since they contain both higher-order self-interaction effects and coupling-induced corrections.

\section{Asymptotic modelling}\label{sec:Asymptotic modelling}

The resonance condition introduced in Section~\ref{sec:Mathematical setting} depends on the spectral properties of the volume integral operator. In the subwavelength regime, these spectral quantities can be approximated by asymptotic expansions or by reduced-order models. Such approximations are computationally efficient and provide explicit insight into the dependence of the resonances on the physical and geometric parameters. However, when inserted into the nonlinear resonance condition, even small errors in the approximation of the operator spectrum may lead to visible errors in the predicted resonant frequencies.

The aim of this section is to recall the asymptotic structure that will later guide the learning procedure. In particular, the logarithmic scaling of the two-dimensional Green's function leads to characteristic terms involving $\delta$, $k_0$, and $\log(\delta k_0)$. These terms will play a central role in the construction of physics-informed features and in the interpretation of the symbolic regression formulas.

\subsection{Integral operator and subwavelength scaling}

Let $D\subset \mathbb{R}^2$ denote the resonator domain. The volume integral operator associated with $D$ is (\ref{eq:defKD}).
We consider the regime in which the resonator is a small rescaling of a fixed reference domain $B$:
\begin{equation}
    D=\delta B,
    \qquad
    0<\delta\ll 1.
\end{equation}
Introducing the rescaled variables
$x=\delta \widetilde{x}$,
$y=\delta \widetilde{y}$,
and writing $\widetilde{u}(\widetilde{y})=u(\delta \widetilde{y})$, we obtain
\begin{equation}\label{eq:operator_rescaling}
    \mathcal{K}_D^\omega[u](\delta \widetilde{x})
    =
    \delta^2
    \int_B
    \left(
        -\frac{1}{2\pi}\log\delta
        -\frac{1}{2\pi}\log|\widetilde{x}-\widetilde{y}|
        +
        R(\delta\widetilde{x},\delta\widetilde{y};\omega)
    \right)
    \widetilde{u}(\widetilde{y})\upd\widetilde{y}.
\end{equation}
This formula separates the explicit dependence on the size parameter $\delta$ from the operator acting on the fixed domain $B$. In particular, it shows that the two-dimensional subwavelength regime is governed not only by powers of $\delta$, but also by logarithmic factors. Moreover, since the regular part of the Green's function depends on the frequency through the background wavenumber $k_0=\omega\sqrt{\varepsilon_0\mu_0}$, higher-order terms involve both $\delta$ and $k_0$.

\subsection{Asymptotic structure of the spectrum}

Let $\lambda(\delta,\omega)$ denote an eigenvalue of $\mathcal{K}_D^\omega$. In the subwavelength regime, the leading spectral behaviour is determined by the logarithmic kernel on the reference domain. The corresponding static logarithmic operator is
\begin{equation}
    \mathcal{K}_B^{(0)}[u](x)
    =
    -\frac{1}{2\pi}
    \int_B \log|x-y|\,u(y)\upd y.
\end{equation}
It appears naturally in the rescaled expression \eqref{eq:operator_rescaling}. The remaining terms in the Green's function produce higher-order corrections that depend on the frequency and on the background wavenumber.

Schematically, the relevant eigenvalues may be written as
\begin{equation}\label{eq:eigenvalue_asymptotic_structure}
    \lambda(\delta,\omega)
    =
    \lambda_{\mathrm{lead}}(\delta)
    +
    \lambda_{\mathrm{corr}}(\delta,\omega),
\end{equation}
where $\lambda_{\mathrm{lead}}$ is determined by the static logarithmic operator and the geometry of $B$, while $\lambda_{\mathrm{corr}}$ contains the frequency-dependent higher-order contributions.

The next-order structure in two dimensions contains terms of the form
\begin{equation}
\label{eq:next_order_structure}
    \delta^2 k_0^2\log(\delta k_0).
\end{equation}
At the operator level, this corresponds schematically to an expansion of the rescaled operator of the form
\begin{equation}\label{eq:higher_order_operator}
    \mathcal{K}_{B}^{\omega}
    \simeq
    \Big(\log(\delta k_0) + \gamma_{\mathrm{E}}\Big)\,\mathcal{K}_{B}^{(-1)}
    +
    \mathcal{K}_{B}^{(0)}
    +
    \delta^2 k_0^2\log(\delta k_0)\,\mathcal{K}_{B}^{(1)}
    + \cdots,
\end{equation}
where \(\mathcal{K}_{B}^{(-1)}\), \(\mathcal{K}_{B}^{(0)}\), and \(\mathcal{K}_{B}^{(1)}\) are operators on the reference domain \(B\), and \(\gamma_{\mathrm{E}}\) denotes the Euler--Mascheroni constant. The precise form of these operators depends on the expansion of the Green's function, but the important point for the present work is the scale of the correction. We refer to \cite{AD} for the expressions of these operators.

In the continuous resonance problem, $k_0$ denotes the background wavenumber and is generally complex when evaluated at a complex resonance. In the learning features, however, we use a real positive frequency scale extracted from the asymptotic resonance. For the single-resonator case, we define
\begin{equation}\label{eq:q_single_def}
    q_{\mathrm{asymp}}
    =
    |\omega_{\mathrm{asymp}}|.
\end{equation}
For the dimer, the corresponding branch-dependent quantities are
\begin{equation}\label{eq:q_dimer_def}
    q_j
    =
    |\omega_{j,\mathrm{asymp}}|,
    \qquad j=1,2.
\end{equation}
Thus, the logarithmic frequency-dependent features used involve real logarithms of positive quantities. In particular, the feature-level analogues of the next-order scales are
\begin{equation}
    \delta^2 q_{\mathrm{asymp}}^2\log (\delta q_{\mathrm{asymp}}),
    \qquad
    \delta^2 q_{\mathrm{asymp}}^2,
\end{equation}
and similarly with $q_j$ for the dimer branches.

Using
\begin{equation}\label{eq:scale}
    \delta^2 q_{\mathrm{asymp}}^2\log(\delta q_{\mathrm{asymp}})
    =
    \delta^2 q_{\mathrm{asymp}}^2\log\delta
    +
    \delta^2 q_{\mathrm{asymp}}^2\log q_{\mathrm{asymp}},
\end{equation}
we obtain the independent feature family
\begin{equation}\label{eq:next_order_feature_family}
    \delta^2 q_{\mathrm{asymp}}^2,
    \qquad
    \delta^2 q_{\mathrm{asymp}}^2\log\delta,
    \qquad
    \delta^2 q_{\mathrm{asymp}}^2\log q_{\mathrm{asymp}}.
\end{equation}
These quantities will be used in the learning framework as physics-informed features and in the assisted symbolic-regression experiments. The combined term $\delta^2 q_{\mathrm{asymp}}^2\log(\delta q_{\mathrm{asymp}})$ is useful for identifying the asymptotic scale, but it is not treated as an independent feature because of the decomposition \eqref{eq:scale}.

\subsection{From operator corrections to resonance residuals}

The appearance of the same feature family in the frequency residual can be justified from the resonance condition. For a given branch, define the scalar resonance function
\begin{equation}
    \Phi(\omega,\delta)
    =
    1+\omega^2\xi(\omega)\lambda(\delta,\omega),
\end{equation}
where $\lambda(\delta,\omega)$ denotes the corresponding eigenvalue of the volume integral operator. The asymptotic resonance $\omega_{\mathrm{asymp}}$ is obtained by replacing $\lambda$ with its asymptotic approximation $\lambda_{\mathrm{asymp}}$. Thus
\begin{equation}
    1+\omega_{\mathrm{asymp}}^2
    \xi(\omega_{\mathrm{asymp}})
    \lambda_{\mathrm{asymp}}(\delta,\omega_{\mathrm{asymp}})
    =0.
\end{equation}
Writing
\begin{equation}
    \lambda(\delta,\omega)
    =
    \lambda_{\mathrm{asymp}}(\delta,\omega)
    +
    \Delta\lambda(\delta,\omega),
\end{equation}
and assuming that the resonance is simple, so that
\begin{equation}
    \partial_\omega \Phi_{\mathrm{asymp}}
    (\omega_{\mathrm{asymp}},\delta)
    \neq 0,
\end{equation}
the implicit function theorem gives, to first order,
\begin{equation}
    \Delta\omega
    =
    \omega_{\mathrm{ref}}-\omega_{\mathrm{asymp}}
    \simeq
    -
    \frac{
    \omega_{\mathrm{asymp}}^2
    \xi(\omega_{\mathrm{asymp}})
    \Delta\lambda(\delta,\omega_{\mathrm{asymp}})
    }
    {
    \partial_\omega \Phi_{\mathrm{asymp}}
    (\omega_{\mathrm{asymp}},\delta)
    }.
\end{equation}
Consequently, the frequency residual inherits the scale of the neglected terms in the eigenvalue and operator expansions, up to multiplication by a smooth branch-dependent factor. Since the next-order operator correction contains the two-dimensional frequency-dependent scale
\begin{equation}
    \delta^2 q_{\mathrm{asymp}}^2\log(\delta q_{\mathrm{asymp}})
    =
    \delta^2 q_{\mathrm{asymp}}^2\log\delta
    +
    \delta^2 q_{\mathrm{asymp}}^2\log q_{\mathrm{asymp}},
\end{equation}
it is natural to include the feature family
\begin{equation}
    \delta^2 q_{\mathrm{asymp}}^2,
    \qquad
    \delta^2 q_{\mathrm{asymp}}^2\log\delta,
    \qquad
    \delta^2 q_{\mathrm{asymp}}^2\log q_{\mathrm{asymp}}
\end{equation}
in the residual model.

\subsection{Asymptotic resonance approximation}

Substituting an asymptotic approximation of $\lambda(\delta,\omega)$ into the resonance condition
$
    1+\omega^2\xi(\omega)\lambda(\delta,\omega)=0
$
yields an approximation of the resonant frequency. At leading order, one obtains a reduced nonlinear equation of the form
$
    1+\omega^2\xi(\omega)\lambda_{\mathrm{asymp}}(\delta,\omega)=0
$,
where $\lambda_{\mathrm{asymp}}$ denotes the chosen asymptotic approximation of the relevant eigenvalue. We denote the resulting resonance by
\begin{equation}
    \omega_{\mathrm{asymp}}.
\end{equation}
The reference resonance, computed from the more accurate spectral model, is denoted by $\omega_{\mathrm{ref}}$.

The discrepancy
\begin{equation}
    \Delta\omega
    =
    \omega_{\mathrm{ref}}-\omega_{\mathrm{asymp}}
\end{equation}
contains the effect of neglected higher-order terms in the operator expansion, as well as the nonlinear propagation of these errors through the resonance condition. The purpose of the data-driven correction developed below is to approximate this residual while preserving the asymptotic structure of the problem.

\subsection{Coupled resonators and interaction operators}

We now consider the dimer configuration $D_{\delta,\kappa}=\delta B_\kappa$, where $B_\kappa$ is the reference dimer made of two unit disks separated by the dimensionless gap $\kappa$. For notational simplicity, we again write $D=D_1\cup D_2$. In this case, the volume integral operator has the block structure (\ref{eq:defKDomegadimer})
where, 
\begin{equation}
    \mathcal{K}_{D_jD_l}^\omega[u](x)
    =
    \int_{D_l}G(x,y;\omega)u(y)\,dy,
    \qquad x\in D_j, \qquad j,l=1,2.
\end{equation}
The diagonal blocks describe the self-interaction of each resonator, while the off-diagonal blocks describe the interaction between the two components.

In the subwavelength regime, the leading behaviour is governed by approximately constant modes on each resonator. This leads to a finite-dimensional reduced model, whose entries are averaged interaction terms. A representative form is
\begin{equation}
    M(\kappa)
    =
    \begin{pmatrix}
        \langle \mathcal{K}_{D_1D_1}^\omega \mathbbm{1}_{D_1},\mathbbm{1}_{D_1}\rangle
        &
        \langle \mathcal{K}_{D_1D_2}^\omega \mathbbm{1}_{D_2},\mathbbm{1}_{D_1}\rangle
        \\[0.2cm]
        \langle \mathcal{K}_{D_2D_1}^\omega \mathbbm{1}_{D_1},\mathbbm{1}_{D_2}\rangle
        &
        \langle \mathcal{K}_{D_2D_2}^\omega \mathbbm{1}_{D_2},\mathbbm{1}_{D_2}\rangle
    \end{pmatrix}.
\end{equation}
The eigenvalues of this reduced matrix approximate the leading eigenvalues of the full operator. The off-diagonal terms depend on the separation distance $\kappa$ and induce a splitting into two branches, associated with symmetric and antisymmetric modes. We denote the corresponding asymptotic resonances by
\begin{equation}
    \omega_{1,\mathrm{asymp}},
    \qquad
    \omega_{2,\mathrm{asymp}}.
\end{equation}

The corresponding reference resonances are denoted by $\omega_{1,\mathrm{ref}}$ and $\omega_{2,\mathrm{ref}}$, and the residuals are
\begin{equation}
    \Delta\omega_j
    =
    \omega_{j,\mathrm{ref}}
    -
    \omega_{j,\mathrm{asymp}},
    \qquad j=1,2.
\end{equation}
Compared with the single-resonator case, these residuals contain both self-interaction corrections and coupling-induced errors. This makes the dimer a useful test case for assessing whether the learned correction can capture interaction effects.

\subsection{Limitations of asymptotic and reduced models}

The asymptotic and reduced models described above provide efficient approximations of the resonances, but they have several limitations. First, their accuracy deteriorates as the system moves away from the strict subwavelength regime. Second, the next terms in the expansion, although known structurally, are difficult to compute explicitly in a form that is convenient for numerical prediction. Third, in coupled systems, interaction effects introduce branch-dependent corrections that are not fully captured by the leading reduced model. Finally, the nonlinear dependence of the resonance condition on $\omega$ can amplify small spectral errors.

These observations motivate a correction strategy that keeps the asymptotic approximation as the leading model, but learns the residual
$
    \omega_{\mathrm{ref}}-\omega_{\mathrm{asymp}}
$
from data. The next section introduces this residual learning framework. The asymptotic structures identified above, in particular the terms involving $\delta^2 q_{\mathrm{asymp}}^2\log(\delta q_{\mathrm{asymp}})$,
will be used to design the input features and to interpret the symbolic formulas obtained later.

\section{Learning the asymptotic residual}\label{sec:Learning framework}

\subsection{Residual formulation}

The asymptotic and reduced models described in Section~\ref{sec:Asymptotic modelling} provide efficient approximations of the resonant frequencies, but they neglect higher-order terms and interaction effects. Rather than learning the resonances directly, we use these approximations as a baseline and learn only the remaining discrepancy.

Let $\omega_{\mathrm{asymp}}$ denote the resonance predicted by an asymptotic or reduced model, and let $\omega_{\mathrm{ref}}$ denote the corresponding reference value. We define the residual
\begin{equation}
    \Delta\omega
    =
    \omega_{\mathrm{ref}}-\omega_{\mathrm{asymp}}.
\end{equation}
The learning task is then to approximate the map
\begin{equation}
    \mathcal{F}:
    X
    \longmapsto
    \Delta\omega,
\end{equation}
where $X$ denotes a vector of physical, asymptotic, and physics-informed features. Once an approximation $\Delta\omega_{\mathrm{ML}}$ has been obtained, the corrected resonance is given by
\begin{equation}
    \omega_{\mathrm{pred}}
    =
    \omega_{\mathrm{asymp}}+\Delta\omega_{\mathrm{ML}}.
\end{equation}

This formulation has two main advantages. First, the leading-order behaviour is already encoded in $\omega_{\mathrm{asymp}}$, so that the learning problem concerns only a smaller correction. Second, the residual is expected to inherit the structure of the neglected asymptotic terms. This makes it possible to use simple learning models, provided that the input features are chosen consistently with the underlying asymptotics.

For the dimer, the same construction is applied branch by branch:
\begin{equation}
    \Delta\omega_j
    =
    \omega_{j,\mathrm{ref}}
    -
    \omega_{j,\mathrm{asymp}},
    \qquad j=1,2.
\end{equation}
The corrected predictions are then
\begin{equation}
    \omega_{j,\mathrm{pred}}
    =
    \omega_{j,\mathrm{asymp}}
    +
    \Delta\omega_{j,\mathrm{ML}},
    \qquad j=1,2.
\end{equation}

\subsection{Feature design}

The choice of features is guided by the asymptotic structure of the problem. We use three types of input quantities.

First, we include the raw physical and geometric parameters, such as the size parameter $\delta$, the damping parameter $\gamma$, and, in the dimer case, the separation distance $\kappa$. Logarithmic transforms of these parameters are also included when they are relevant to the asymptotic scaling.

Second, we include quantities obtained from the asymptotic or reduced model itself. These include the real and imaginary parts of $\omega_{\mathrm{asymp}}$, the approximate eigenvalues $\lambda_{\mathrm{asymp}}$, and, for the dimer, quantities measuring the splitting between the two branches. These features encode the leading-order spectral information already available from the reduced model.

Third, we include physics-informed features motivated by the next terms in the asymptotic expansion. In the single-resonator case, these include terms such as
\begin{equation}
    \delta^2,
    \qquad
    \delta^2\log\delta,
    \qquad
    \gamma\delta^2,
\end{equation}
together with the frequency-dependent family
\begin{equation}
    \delta^2 q_{\mathrm{asymp}}^2,
    \qquad
    \delta^2 q_{\mathrm{asymp}}^2\log\delta,
    \qquad
    \delta^2 q_{\mathrm{asymp}}^2\log q_{\mathrm{asymp}},
\end{equation}
Here $q_{\mathrm{asymp}}$ is defined by \eqref{eq:q_single_def}. For the dimer, the same features are evaluated branchwise using $q_j$ from \eqref{eq:q_dimer_def}. This notation separates the physical background wavenumber $k_0$, which belongs to the resonance problem, from the positive real frequency scale used as an input feature in the learning models.
These features reflect the structure of the two-dimensional Green's function expansion discussed in Section~\ref{sec:Asymptotic modelling}. In the dimer case, analogous branch-dependent features are used, evaluated with the corresponding asymptotic wavenumbers, together with interaction features involving the eigenvalue splitting and the separation distance $\kappa$.

The purpose of this feature design is not to replace the asymptotic model, but to expose the dominant scales of the missing correction to the learning algorithm. This will also allow us to test, through the ablation study, whether the improvement genuinely comes from physically meaningful structure.

\subsection{Regression model}

As a first correction model, we use Ridge regression. The real and imaginary parts of the residual are learned simultaneously. Given a feature vector $X\in\mathbb{R}^p$, we approximate
\begin{equation}
    \Delta\omega
    \approx
    WX+b,
\end{equation}
where $W\in\mathbb{R}^{2\times p}$ and $b\in\mathbb{R}^2$. The two components correspond to
$
    \Delta\omega
    =
    \big(\operatorname{Re}\Delta\omega,
    \operatorname{Im}\Delta\omega\big)
$.
The parameters are obtained by minimizing the regularized least-squares functional
\begin{equation}
    \min_{W,b}
    \sum_{n=1}^N
    \left\|
        WX_n+b-\Delta\omega_n
    \right\|^2
    +
    \alpha_{\mathrm{R}}
    \|W\|^2,
\end{equation}
where $\alpha_{\mathrm{R}}>0$ is the Ridge regularization parameter. The numerical details of the Ridge implementation, including feature standardization, train-test splits, and the value of the regularization parameter, are given in Appendix~\ref{app:numerical_details}.

The use of a linear model is deliberate. Since the asymptotic approximation already captures the dominant behaviour, and since the features include the expected higher-order scalings, the residual should be well approximated by a low-dimensional combination of these quantities. This also makes the learned correction easier to interpret than a high-capacity black-box model.

The model is trained on samples of the form (features, $\omega_{\mathrm{asymp}}$, $\omega_{\mathrm{ref}}$), with target $\omega_{\mathrm{ref}}-\omega_{\mathrm{asymp}}$. Details on the train-test split and preprocessing are given in Appendix~\ref{app:numerical_details}.

\subsection{Corrected prediction and error measure}

Given a trained model, the corrected prediction is
\begin{equation}
    \omega_{\mathrm{pred}} = \omega_{\mathrm{asymp}} + \Delta\omega_{\mathrm{ML}}.
\end{equation}
The accuracy is measured by the relative error
\begin{equation}
    E(\omega) = \frac{|\omega_{\mathrm{pred}}-\omega_{\mathrm{ref}}|}{|\omega_{\mathrm{ref}}|}.
\end{equation}
The same definition is used to evaluate the uncorrected asymptotic model, replacing $\omega_{\mathrm{pred}}$ by $\omega_{\mathrm{asymp}}$. This allows a direct comparison between the original asymptotic approximation and the learned correction.

\subsection{Interpretation}

The residual learning framework can be viewed as a data-driven extension of asymptotic modelling. The asymptotic approximation provides the leading-order physics, while the learned residual captures higher-order terms and coupling effects that are difficult to compute explicitly.

The Ridge model provides a stable and accurate correction, but its output is still a fitted linear combination of many features. In the next section, we go one step further and use symbolic regression to search for compact analytic expressions for the residual. This provides a more interpretable correction and makes it possible to compare the learned formulas directly with the expected asymptotic structure.

\section{Symbolic regression of the residual}\label{sec:symbolic_regression}

The residual learning framework introduced in Section~\ref{sec:Learning framework} improves the accuracy of the asymptotic approximation while preserving the leading-order physical structure. However, even when a simple linear model is used, the learned correction remains a fitted combination of prescribed features. In order to obtain more explicit and interpretable correction formulas, we complement the regression approach with symbolic regression.

Symbolic regression seeks analytic expressions that approximate a target quantity from data by combining input variables through elementary operations. Unlike standard regression models, which fit coefficients within a fixed functional form, symbolic regression searches over the space of mathematical expressions. This makes it particularly useful in scientific applications, where the goal is not only to predict accurately but also to identify compact formulas that reveal structure in the data. In this work, we use PySR, based on the SymbolicRegression.jl framework, following the approach described in \cite{cranmer2023interpretable}.

\subsection{Symbolic regression target}

We apply symbolic regression to the asymptotic residual rather than to the full resonance. For a single resonator, the target is
\begin{equation}
    \Delta\omega
    =
    \omega_{\mathrm{ref}}
    -
    \omega_{\mathrm{asymp}}.
\end{equation}
For the dimer, the same procedure is applied separately to the two branches:
\begin{equation}
    \Delta\omega_j
    =
    \omega_{j,\mathrm{ref}}
    -
    \omega_{j,\mathrm{asymp}},
    \qquad j=1,2.
\end{equation}
Since the resonances are complex-valued, we fit the real and imaginary parts separately. Thus, for each configuration, symbolic regression is used to find expressions for
$\operatorname{Re}\Delta\omega$ and
$\operatorname{Im}\Delta\omega$.
This choice is important. Learning the full resonance directly would require the symbolic model to rediscover the leading-order asymptotic behaviour. By learning only the residual, the symbolic search is focused on the missing higher-order correction. This is consistent with the philosophy of the paper: the asymptotic model provides the dominant physical approximation, while data-driven methods are used only to model the discrepancy.

\subsection{Input variables and asymptotic structure}

The input variables used for symbolic regression are chosen from the same physics-informed feature set as in Section~\ref{sec:Learning framework}. For the single-resonator problem, these include the size parameter $\delta$, the damping parameter $\gamma$, the asymptotic resonance $\omega_{\mathrm{asymp}}$, and asymptotically motivated combinations such as
\begin{equation}
    \delta^2,
    \qquad
    \delta^2 q_{\mathrm{asymp}}^2,
    \qquad
    \delta^2 q_{\mathrm{asymp}}^2\log q_{\mathrm{asymp}},
    \qquad
    \delta^2 q_{\mathrm{asymp}}^2\log(\delta).
\end{equation}
For the dimer, the symbolic regression also uses branch-dependent asymptotic quantities, including the approximate eigenvalues and the splitting between the two branches. If $\lambda_{1,\mathrm{asymp}}$ and $\lambda_{2,\mathrm{asymp}}$ denote the two leading asymptotic eigenvalues, we use in particular the splitting
\begin{equation}
    \Lambda
    =
    \left|
        \lambda_{1,\mathrm{asymp}}
        -
        \lambda_{2,\mathrm{asymp}}
    \right|.
\end{equation}

The inclusion of these variables is guided by the asymptotic expansion. In particular, the term $\delta^2 q_{\mathrm{asymp}}^2\log(\delta q_{\mathrm{asymp}})$ appears naturally in the next-order correction to the two-dimensional volume integral operator \eqref{eq:higher_order_operator}. Therefore, when symbolic regression selects expressions involving this family of terms, the result can be interpreted as being consistent with the expected asymptotic structure. We emphasize, however, that the symbolic formulas obtained below are data-driven fitted corrections; they should not be interpreted as rigorous derivations of the full higher-order asymptotic expansion.

\subsection{Discovered symbolic correction formulas}

We now report compact formulas obtained by symbolic regression. The formulas in this subsection are \emph{data-driven symbolic corrections}: they are fitted correction laws for the residual in the sampled parameter regime, not rigorous asymptotic expansions. Their purpose is to expose which combinations of variables are used by the learned correction and to provide explicit analytic approximations of the residual.

All quantities appearing below are non-dimensional, following the convention introduced in Section~\ref{sec:Mathematical setting}. In particular, $\omega$, $\gamma$, $\lambda_j$, $q_{\mathrm{asymp}}$, and $q_j$ denote non-dimensional quantities, and the numerical constants in the symbolic expressions are dimensionless fitted coefficients. The quantities $q_{\mathrm{asymp}}$ and $q_j$ are the positive real frequency scales defined in \eqref{eq:q_single_def} and \eqref{eq:q_dimer_def}. Thus logarithms involving $q_{\mathrm{asymp}}$ or $q_j$ are real logarithms of positive quantities.

It is important to separate these symbolic formulas from the asymptotic expansion itself. Symbolic regression does not derive the next term in the asymptotic expansion. Rather, it searches for compact data-driven expressions that are compatible with the asymptotically motivated feature scales. Schematically, in the assisted setting one may view the symbolic correction as
\begin{equation}
    \Delta\omega_{\mathrm{SR}}
    =
    \underbrace{
        \Delta\omega_{\mathrm{proxy}}
    }_{\text{motivated by asymptotics}}
    +
    \underbrace{
        \Delta\omega_{\mathrm{data}}
    }_{\text{discovered from data}},
\end{equation}
where $\Delta\omega_{\mathrm{proxy}}$ is built from the feature scales suggested by the operator expansion, while $\Delta\omega_{\mathrm{data}}$ is the remaining symbolic correction fitted from data. In the unassisted symbolic-regression experiments, the whole correction is fitted directly from data; in the assisted experiments, these two contributions are separated explicitly.

\subsubsection*{Single resonator}

For the single-resonator problem, one compact symbolic correction is obtained as follows. Writing
\begin{equation}
    \omega_{\mathrm{asymp}}
    =
    \omega_{\mathrm{asymp},r}
    +
    i\omega_{\mathrm{asymp},i},
\end{equation}
the corrected resonance is
\begin{equation}
    \omega_{\mathrm{sym}}
    =
    \omega_{\mathrm{asymp}}
    +
    \operatorname{Re}(\Delta\omega_{\mathrm{sym}})
    +
    i\operatorname{Im}(\Delta\omega_{\mathrm{sym}}),
\end{equation}
where
\begin{equation}\label{eq:single_symbolic_real}
    \operatorname{Re}(\Delta\omega_{\mathrm{sym}})
    =
    \delta^2\omega_{\mathrm{asymp},r}
    \left(
        0.02398765\,\omega_{\mathrm{asymp},r}^2
        +
        0.00012753597\,\omega_{\mathrm{asymp},r}\omega_{\mathrm{asymp},i}
        +
        0.0053167343
    \right),
\end{equation}
and
\begin{equation}\label{eq:single_symbolic_imag}
    \operatorname{Im}(\Delta\omega_{\mathrm{sym}})
    =
    0.02947997\,\omega_{\mathrm{asymp},i}
    \left(
        \delta^2 q_{\mathrm{asymp}}^2
        +
        \delta^2
        +
        \delta^2 q_{\mathrm{asymp}}^2\log q_{\mathrm{asymp}}
    \right).
\end{equation}

The real part is dominated by a $\delta^2$ correction modulated by the leading asymptotic resonance. The imaginary part contains the frequency-dependent structures $\delta^2 q_{\mathrm{asymp}}^2$ and $\delta^2 q_{\mathrm{asymp}}^2\log q_{\mathrm{asymp}}$, which are consistent with the feature scales motivated by the next-order operator expansion. The factor $\omega_{\mathrm{asymp},i}$ also reflects that the correction to the imaginary part is tied to the damping of the resonance.

\subsubsection*{Dimer}

For the dimer, symbolic regression is applied separately to the two resonance branches. Let
\begin{equation}
    \Lambda
    =
    \left|
        \lambda_{1,\mathrm{asymp}}
        -
        \lambda_{2,\mathrm{asymp}}
    \right|
\end{equation}
denote the non-dimensional asymptotic splitting of the two leading eigenvalues.

For the first branch, the corrected resonance is
\begin{equation}
    \omega_{1,\mathrm{sym}}
    =
    \omega_{1,\mathrm{asymp}}
    +
    \operatorname{Re}(\Delta\omega_{1,\mathrm{sym}})
    +
    i\operatorname{Im}(\Delta\omega_{1,\mathrm{sym}}),
\end{equation}
with
\begin{equation}\label{eq:dimer_branch1_symbolic_real}
    \operatorname{Re}(\Delta\omega_{1,\mathrm{sym}})
    =
    -0.0298\,\delta^2 q_1^2
    \left(
        \Lambda
        -
        0.285
    \right),
\end{equation}
and
\begin{equation}\label{eq:dimer_branch1_symbolic_imag}
    \operatorname{Im}(\Delta\omega_{1,\mathrm{sym}})
    =
    0.034346\,\gamma\delta^2
    \left(
        \Lambda
        -
        0.28841916
    \right).
\end{equation}

For the second branch, the corrected resonance is
\begin{equation}
    \omega_{2,\mathrm{sym}}
    =
    \omega_{2,\mathrm{asymp}}
    +
    \operatorname{Re}(\Delta\omega_{2,\mathrm{sym}})
    +
    i\operatorname{Im}(\Delta\omega_{2,\mathrm{sym}}),
\end{equation}
where
\begin{equation}\label{eq:dimer_branch2_symbolic_real}
    \operatorname{Re}(\Delta\omega_{2,\mathrm{sym}})
    =
    \delta^2 q_2^2
    \left(
        0.061492536
        -
        0.33402312\,\lambda_{2,\mathrm{asymp}}
    \right)
    -
    1.45\times 10^{-6},
\end{equation}
and
\begin{equation}\label{eq:dimer_branch2_symbolic_imag}
    \operatorname{Im}(\Delta\omega_{2,\mathrm{sym}})
    =
    \gamma\delta^2\lambda_{2,\mathrm{asymp}}^2
    \left(
        \operatorname{Re}(\omega_{2,\mathrm{asymp}})
        -
        \lambda_{1,\mathrm{asymp}}
        -
        0.26824647
    \right).
\end{equation}

Since all quantities are non-dimensional, combinations such as
$
    \operatorname{Re}(\omega_{2,\mathrm{asymp}})
    -
    \lambda_{1,\mathrm{asymp}}
$
should be read as combinations of non-dimensional input features selected by symbolic regression.

These expressions reveal two important features. First, the corrections are organized around the expected powers of the asymptotic parameter, especially $\delta^2$ and $\delta^2 q_j^2$. Second, the dimer corrections depend explicitly on branch-dependent spectral quantities, through $\lambda_{j,\mathrm{asymp}}$ and the splitting $\Lambda$. This provides numerical evidence that the learned residual contains interaction-dependent structure, rather than only single-resonator self-interaction corrections.

\subsection{Asymptotically assisted symbolic regression}

We also consider two variants in which symbolic regression is explicitly guided by the next-order asymptotic structure. The purpose of these experiments is not to claim that symbolic regression derives the next asymptotic term, but to test whether supplying this structure makes the data-driven correction simpler and more stable.

We first construct a calibrated next-order proxy, denoted by $\Delta\omega_{\mathrm{proxy}}$, using the feature family
\begin{equation}
    \delta^2 q_{\mathrm{asymp}}^2,
    \qquad
    \delta^2 q_{\mathrm{asymp}}^2\log\delta,
    \qquad
    \delta^2 q_{\mathrm{asymp}}^2\log q_{\mathrm{asymp}},
\end{equation}
and, in the dimer case, the corresponding branch-dependent quantities obtained by replacing $q_{\mathrm{asymp}}$ with $q_j$.

We then compare two assisted strategies.

The first strategy is \emph{asymptotic preconditioning}. In this case, the proxy is subtracted before applying symbolic regression. The symbolic model is trained on
\begin{equation}
    \omega_{\mathrm{ref}}
    -
    \left(
        \omega_{\mathrm{asymp}}
        +
        \Delta\omega_{\mathrm{proxy}}
    \right).
\end{equation}
Thus symbolic regression is asked to learn only the residual left after the asymptotically motivated proxy has been removed.

The second strategy is \emph{asymptotic feature augmentation}. In this case, the symbolic model is still trained on the full residual
\begin{equation}
    \omega_{\mathrm{ref}}-\omega_{\mathrm{asymp}},
\end{equation}
but $\Delta\omega_{\mathrm{proxy}}$ is supplied as an additional input feature. This allows symbolic regression to decide how to use the proxy, rather than forcing the proxy to be subtracted beforehand.

These two strategies test different ways of injecting asymptotic information into the symbolic search. Asymptotic preconditioning changes the target, while asymptotic feature augmentation changes the input representation.
The numerical results in Section~\ref{sec:Numerical experiments} show that this assisted symbolic regression improves the accuracy further, for both the single resonator and the dimer. Thus, symbolic regression is most effective when used not as a purely black-box discovery tool, but as a mechanism for extracting compact formulas from an asymptotically informed feature space.

\section{Numerical experiments}\label{sec:Numerical experiments}

In this section, we evaluate the proposed asymptotics-guided correction framework on single-resonator and dimer configurations. The goal is twofold. First, we quantify the improvement obtained by learning the asymptotic residual with Ridge regression. Second, we assess whether symbolic regression can produce compact formulas that retain the accuracy of the learned correction while improving interpretability.

The reference resonances $\omega_{\mathrm{ref}}$ are computed from the accurate nonlinear spectral model, while $\omega_{\mathrm{asymp}}$ denotes the asymptotic or reduced approximation. Given a learned residual $\Delta\omega_{\mathrm{ML}}$, the corrected prediction is
\begin{equation}
    \omega_{\mathrm{pred}} =
    \omega_{\mathrm{asymp}} + \Delta\omega_{\mathrm{ML}}.
\end{equation}
For symbolic regression, we similarly write
\begin{equation}
    \omega_{\mathrm{sym}} = 
    \omega_{\mathrm{asymp}} + \Delta\omega_{\mathrm{sym}}.
\end{equation}
Accuracy is measured using the relative error
\begin{equation}
    E(\omega) =
    \frac{|\omega_{\mathrm{pred}}-\omega_{\mathrm{ref}}|}{|\omega_{\mathrm{ref}}|},
\end{equation}
with the obvious replacement of $\omega_{\mathrm{pred}}$ by $\omega_{\mathrm{asymp}}$ or $\omega_{\mathrm{sym}}$ when evaluating the asymptotic or symbolic models.

\subsection{Learning protocol and validation strategy}

For each configuration, the dataset is split into training and test sets. The in-distribution test corresponds to a random held-out subset of the sampled parameter domain. To assess extrapolation, we also consider out-of-distribution tests in which the model is trained away from the most challenging regimes and evaluated on held-out parameter ranges. For the single resonator, we test extrapolation toward larger values of $\delta$ and larger damping $\gamma$. For the dimer, we test extrapolation toward larger $\delta$ and smaller separation distance $\kappa$, where interaction effects are stronger.

All pre-processing steps, including feature standardization, are fitted only on the training set and then applied to the test set. The learning target is always the residual $\omega_{\mathrm{ref}}-\omega_{\mathrm{asymp}}$, not the full resonance. We retain only the physical resonance branch satisfying $\operatorname{Re}\omega>0$ and $\operatorname{Im}\omega<0$, which correspond to the outgoing damped resonances considered in this work. In the dimer case, both resonance branches are included in the error averages.

The mean relative error reported in the tables is the average of the relative error 
over the held-out test set:
\begin{align}\label{eq:to be said}
    \overline{E}
    =
    \frac{1}{N_{\mathrm{test}}}
    \sum_{j=1}^{N_{\mathrm{test}}} E_j, \qquad E_j
    =
    \frac{|\omega_{\mathrm{pred},j}-\omega_{\mathrm{ref},j}|}
    {|\omega_{\mathrm{ref},j}|}.
\end{align}
For the dimer, the average is taken over both resonance branches.

\begin{table}[ht]
    \centering
    \begin{tabular}{lccc}
        \hline
        Validation case & Asymptotic mean error & Corrected mean error & Improvement factor \\
        \hline
        Single, in-distribution & $9.09\times 10^{-4}$ & $2.22\times 10^{-7}$ & $4.10\times 10^{3}$ \\
        Single, large-$\delta$ OOD & $2.29\times 10^{-3}$ & $7.89\times 10^{-7}$ & $2.90\times 10^{3}$ \\
        Single, large-$\gamma$ OOD & $9.94\times 10^{-4}$ & $7.92\times 10^{-7}$ & $1.26\times 10^{3}$ \\
        Dimer, in-distribution & $1.22\times 10^{-4}$ & $3.47\times 10^{-6}$ & $35.2$ \\
        Dimer, large-$\delta$ OOD & $3.68\times 10^{-4}$ & $7.35\times 10^{-5}$ & $5.01$ \\
        Dimer, small-$\kappa$ OOD & $1.47\times 10^{-4}$ & $3.83\times 10^{-5}$ & $3.82$ \\
        \hline
    \end{tabular}
    \caption{In-distribution and out-of-distribution validation of the physics-informed residual correction. OOD denotes out-of-distribution. The corrected error corresponds to the Ridge residual model using physics-informed features.}
    \label{tab:in_out_distribution_validation}
\end{table}

The validation table reports the full physics-informed Ridge correction for each validation setting. Table~\ref{tab:symbolic_summary} uses the same in-distribution held-out test sets as the first and fourth rows of Table~\ref{tab:in_out_distribution_validation}, and compares the full physics-informed Ridge correction with symbolic and assisted symbolic regression.

\begin{table}[ht]
    \centering
    \begin{tabular}{lcccc}
        \hline
        Configuration & Asymptotic & Ridge & Symbolic & Assisted symbolic \\
        \hline
        Single resonator & $9.09\times10^{-4}$ & $2.22\times10^{-7}$ & $4.70\times10^{-8}$ & $1.89\times10^{-11}$ \\
        Dimer & $1.22\times10^{-4}$ & $3.47\times10^{-6}$ & $8.35\times10^{-7}$ & $9.30\times10^{-11}$ \\
        \hline
    \end{tabular}
    \caption{Mean relative errors on the in-distribution held-out test sets for the asymptotic approximation, the full physics-informed Ridge correction, symbolic regression, and assisted symbolic regression.}
    \label{tab:symbolic_summary}
\end{table}
These values are obtained on independent held-out test sets. They show that the learned residual is highly structured: Ridge regression already gives a substantial correction, while symbolic and assisted symbolic regression recover compact formulas with significantly smaller errors.

\subsection{Single resonator}

We first consider a single resonator over the parameter ranges described in Appendix~\ref{app:numerical_details}. The baseline approximation is the asymptotic resonance $\omega_{\mathrm{asymp}}$, and the learning target is the residual
$
\Delta\omega = \omega_{\mathrm{ref}} - \omega_{\mathrm{asymp}}
$.

\begin{figure}[ht]
    \centering
    \begin{subfigure}[t]{0.47\textwidth}
        \centering
        \includegraphics[width=\linewidth]{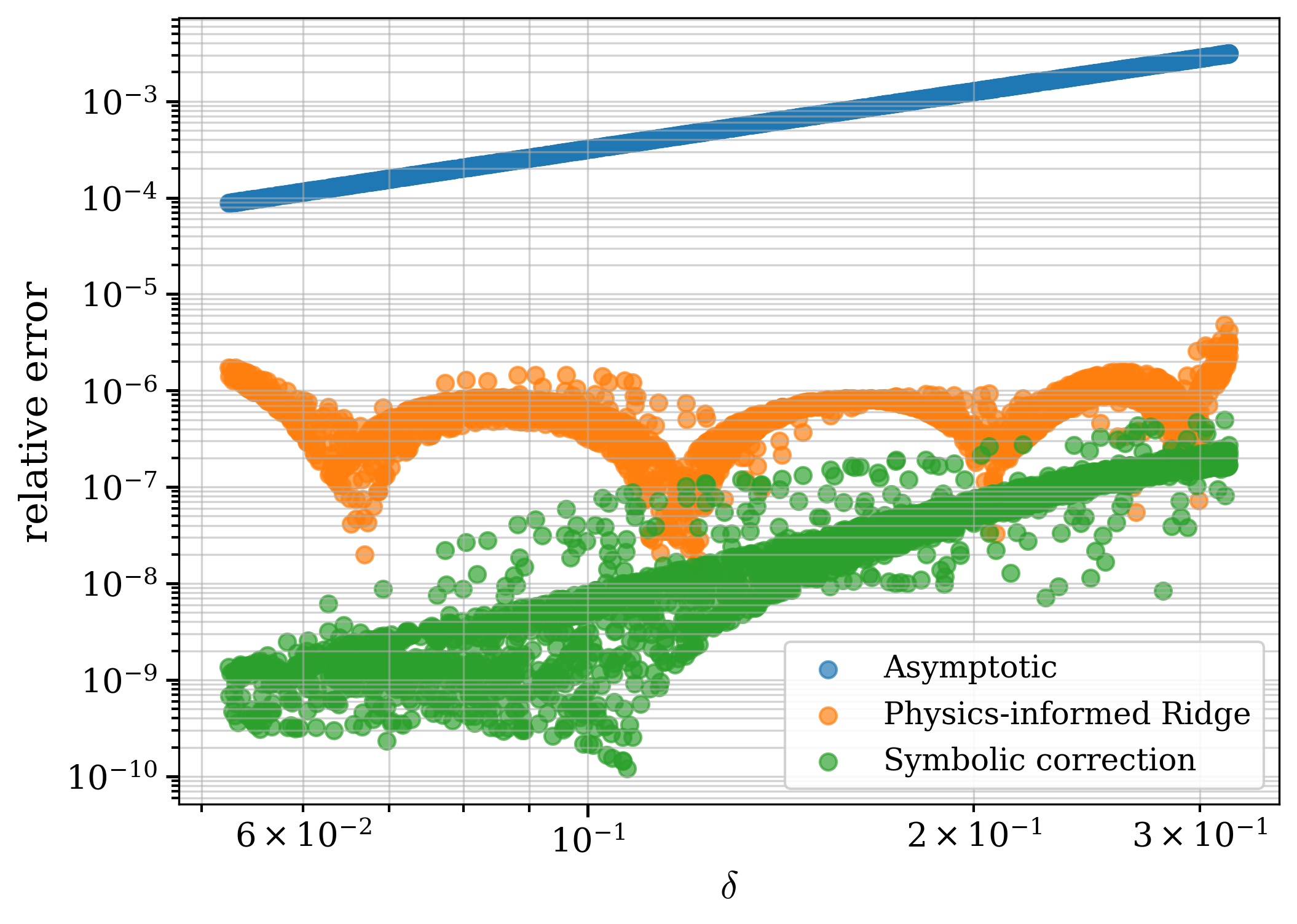}
        \caption{Relative error versus the size parameter $\delta$.}
        \label{fig:single_error_delta}
    \end{subfigure}
    \hfill
    \begin{subfigure}[t]{0.42\textwidth}
        \centering
        \includegraphics[width=\linewidth]{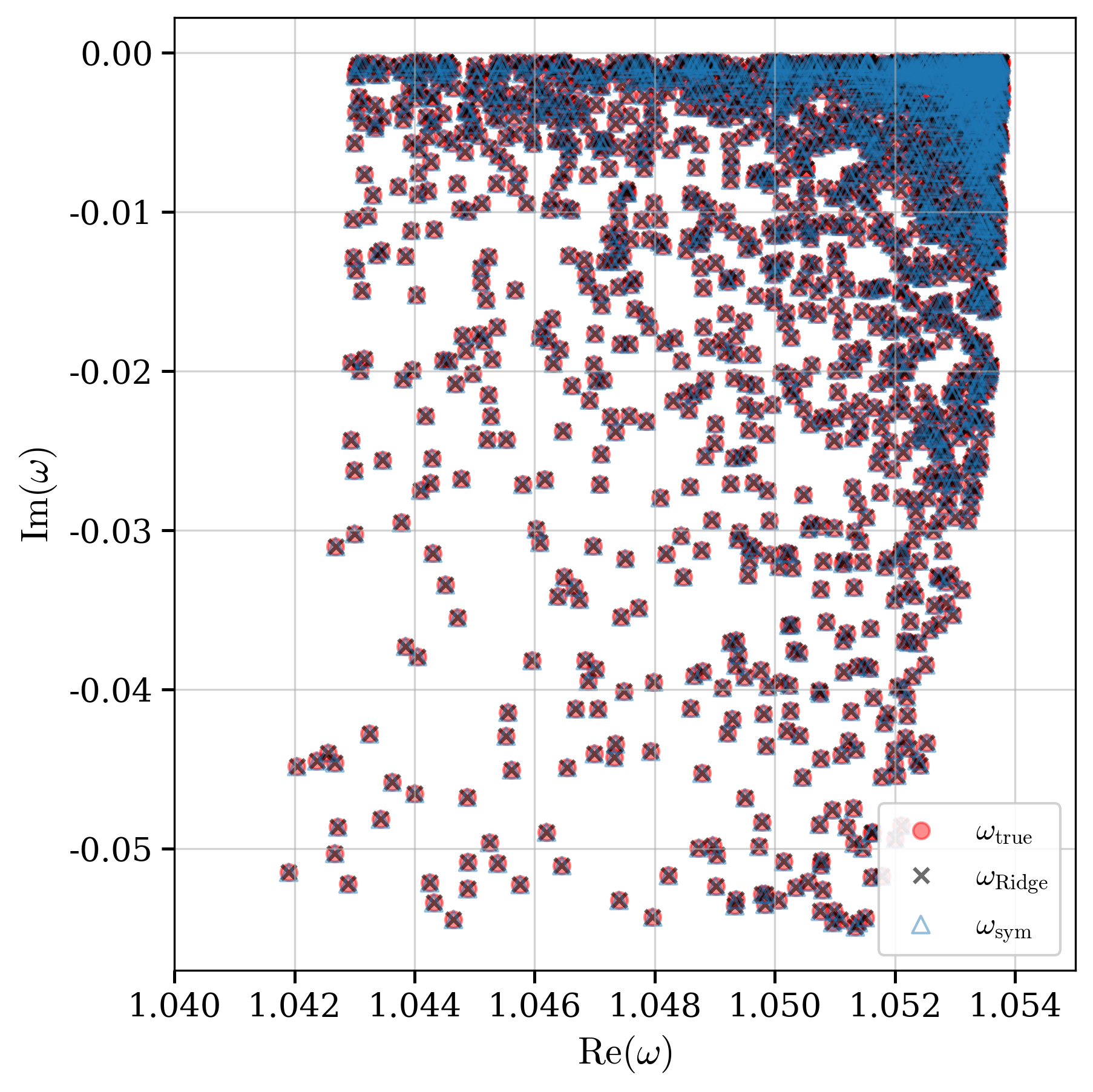}
        \caption{Reconstruction of the resonances in the complex plane.}
        \label{fig:single_complex}
    \end{subfigure}
    \caption{Single-resonator prediction results. Left: relative error of the asymptotic approximation and of the corrected models as a function of the asymptotic parameter $\delta$. Right: comparison of the predicted and reference resonances in the complex plane. The asymptotic approximation deteriorates as $\delta$ increases, while the learned corrections remain close to the reference resonances over the sampled regime.}
    \label{fig:single:error_vs_delta}
\end{figure}

Figure~\ref{fig:single_error_delta} shows the relative error as a function of $\delta$ for the asymptotic approximation and the Ridge-corrected model. As expected, the asymptotic approximation deteriorates as $\delta$ increases, reflecting the increasing importance of higher-order terms. The learned correction significantly reduces the error across the full range of $\delta$, showing that the residual is effectively captured by the physics-informed features. 
The symbolic model further improves the prediction while retaining an explicit analytic form for the residual. In this case, the mean relative error is reduced from approximately $9.09\times 10^{-4}$ for the asymptotic approximation to approximately $4.70\times 10^{-8}$ for the symbolic correction. In Figure~\ref{fig:single_complex}, the reconstruction in the complex plane is shown. The corrected predictions are nearly indistinguishable from the reference resonances, indicating that both the real and imaginary parts are accurately recovered.

\subsection{Dimer configuration}

The dimer configuration is not intended as a simple repetition of the single-resonator experiment. It is the first nontrivial test of whether residual learning transfers from isolated self-interaction corrections to interaction-induced spectral corrections. In this case, the reduced model must approximate not only the resonance of each component, but also the coupling between the two resonators and the resulting branch splitting. The learned residual is therefore branch-dependent and depends on the separation parameter $\kappa$, in addition to $\delta$ and $\gamma$.

The convergence of the dimer reference discretization is checked in Appendix~\ref{app:dimer_convergence}, where the reference computation is compared with finer and coarser discretizations.

\begin{figure}[ht]
    \centering
    \begin{subfigure}[t]{0.43\textwidth}
        \centering
        \includegraphics[width=\linewidth]{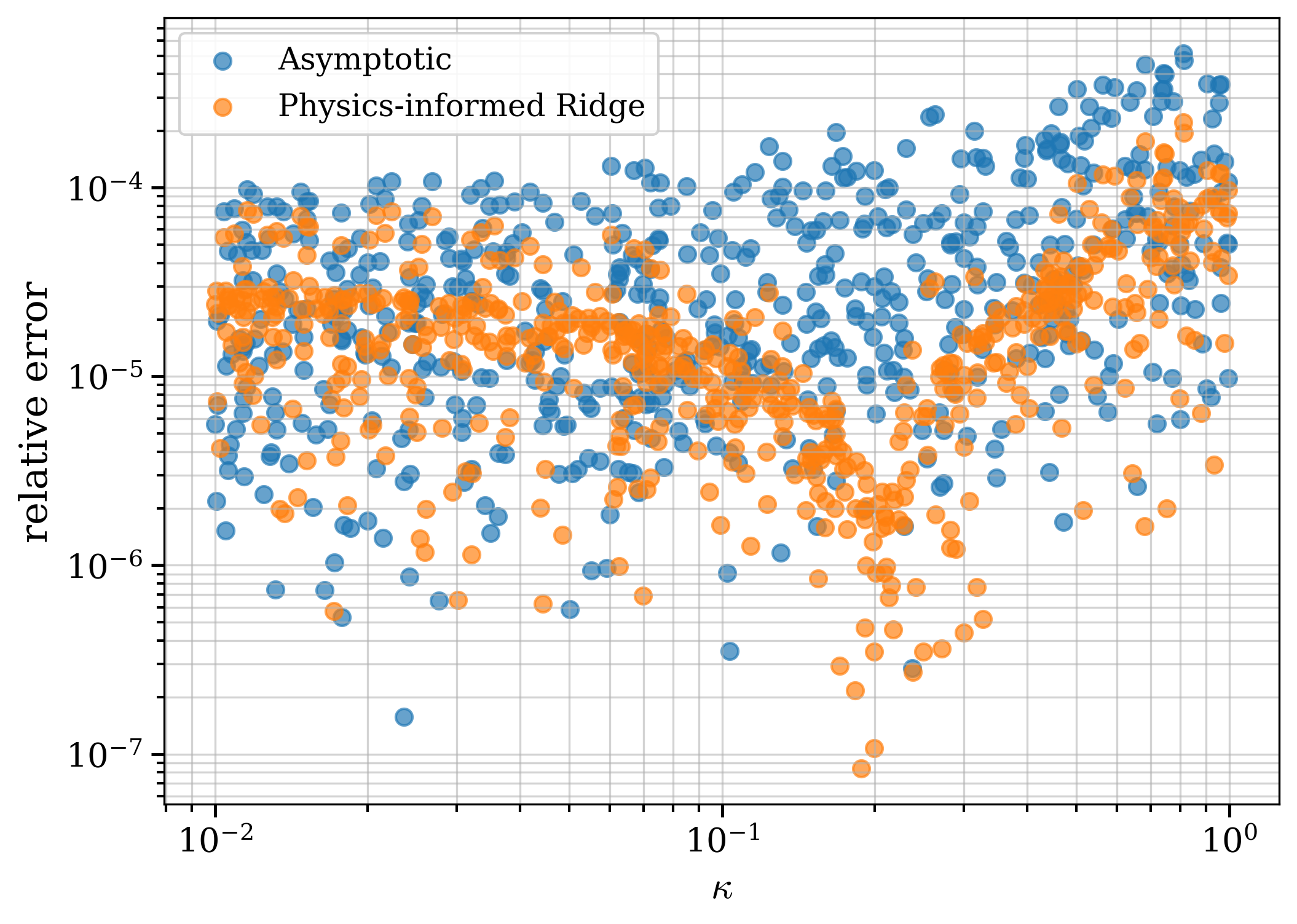}
        \caption{First dimer branch: relative error as a function of the separation distance $\kappa$.}
        \label{fig:dimer_branch1}
    \end{subfigure}
    \hfill
    \begin{subfigure}[t]{0.43\textwidth}
        \centering
        \includegraphics[width=\linewidth]{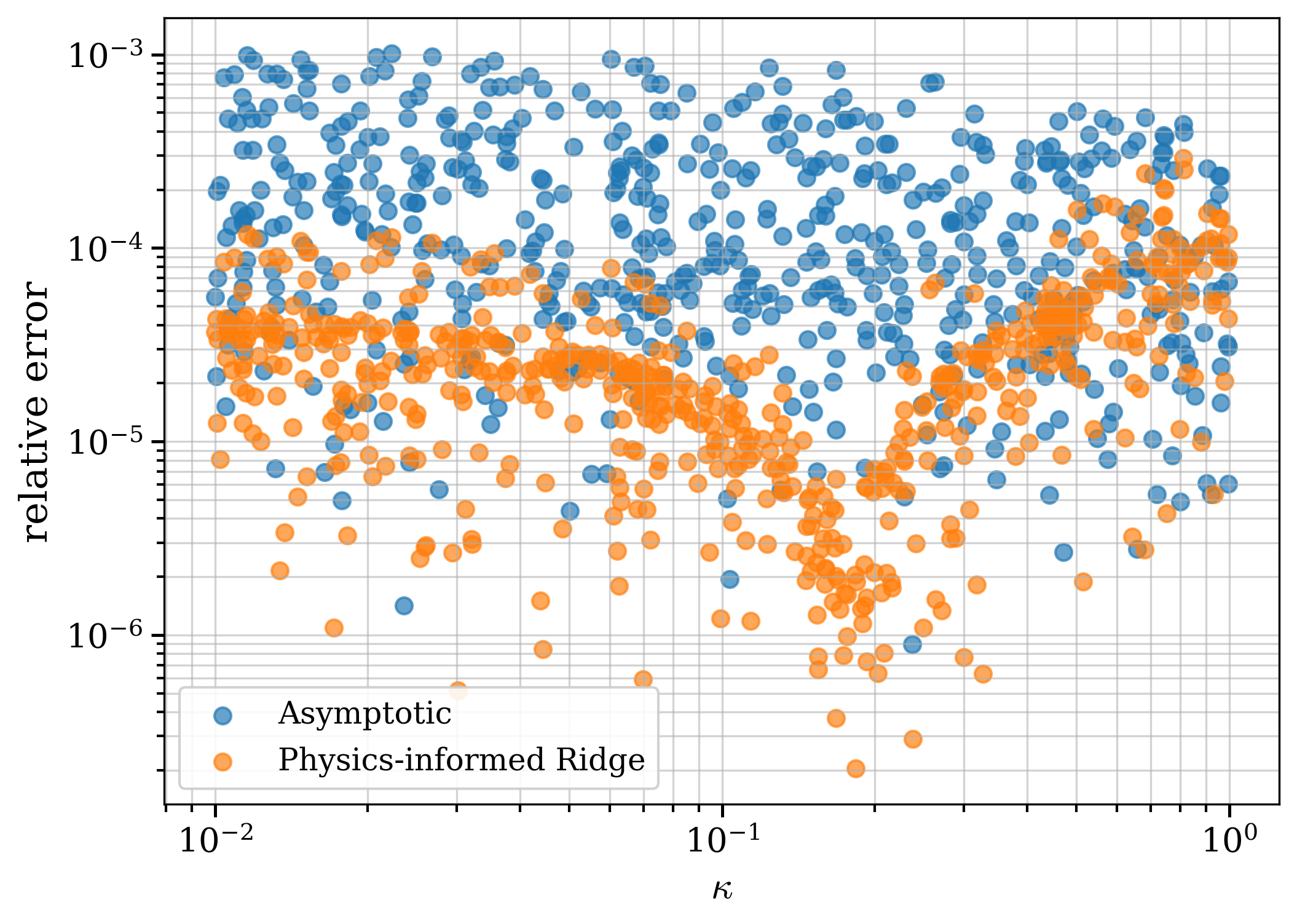}
        \caption{Second dimer branch: relative error as a function of the separation distance $\kappa$.}
        \label{fig:dimer_branch2}
    \end{subfigure}
    \caption{Dimer prediction errors as functions of the separation distance $\kappa$ for the two resonance branches. The reduced asymptotic model captures the leading coupling effect, but the Ridge correction significantly improves the approximation for both branches. The second branch exhibits a larger variability, reflecting its stronger sensitivity to interaction effects.}
    \label{fig:dimer:coarse-to-accurate}
\end{figure}

Figure~\ref{fig:dimer_branch1} shows the relative error for the first branch. The reduced asymptotic model captures the leading interaction effect, but it does not fully account for the branch-dependent correction. The Ridge model substantially reduces the error over the range of separations. Then, Figure~\ref{fig:dimer_branch2} shows the corresponding result for the second branch. This branch is more sensitive to coupling effects, and the baseline error is more variable. Nevertheless, the learned correction remains effective and consistently improves the reduced approximation.

\begin{figure}[ht]
    \centering
    \begin{subfigure}[t]{0.43\textwidth}
        \centering
        \includegraphics[width=\linewidth]{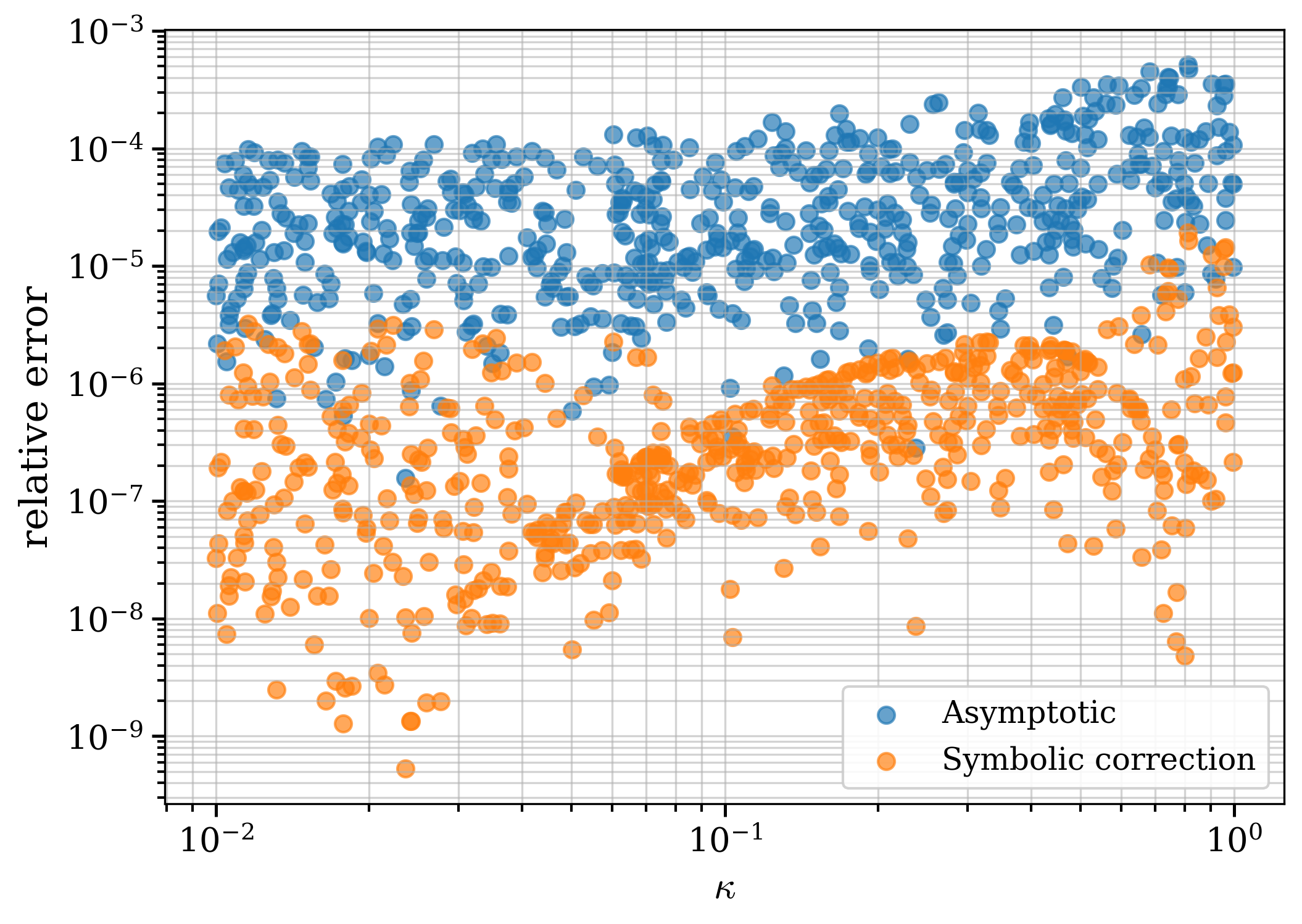}
        \caption{First dimer branch: error of the compact symbolic-regression correction.}
        \label{fig:dimer_symbolic_branch1}
    \end{subfigure}
    \hfill
    \begin{subfigure}[t]{0.43\textwidth}
        \centering
        \includegraphics[width=\linewidth]{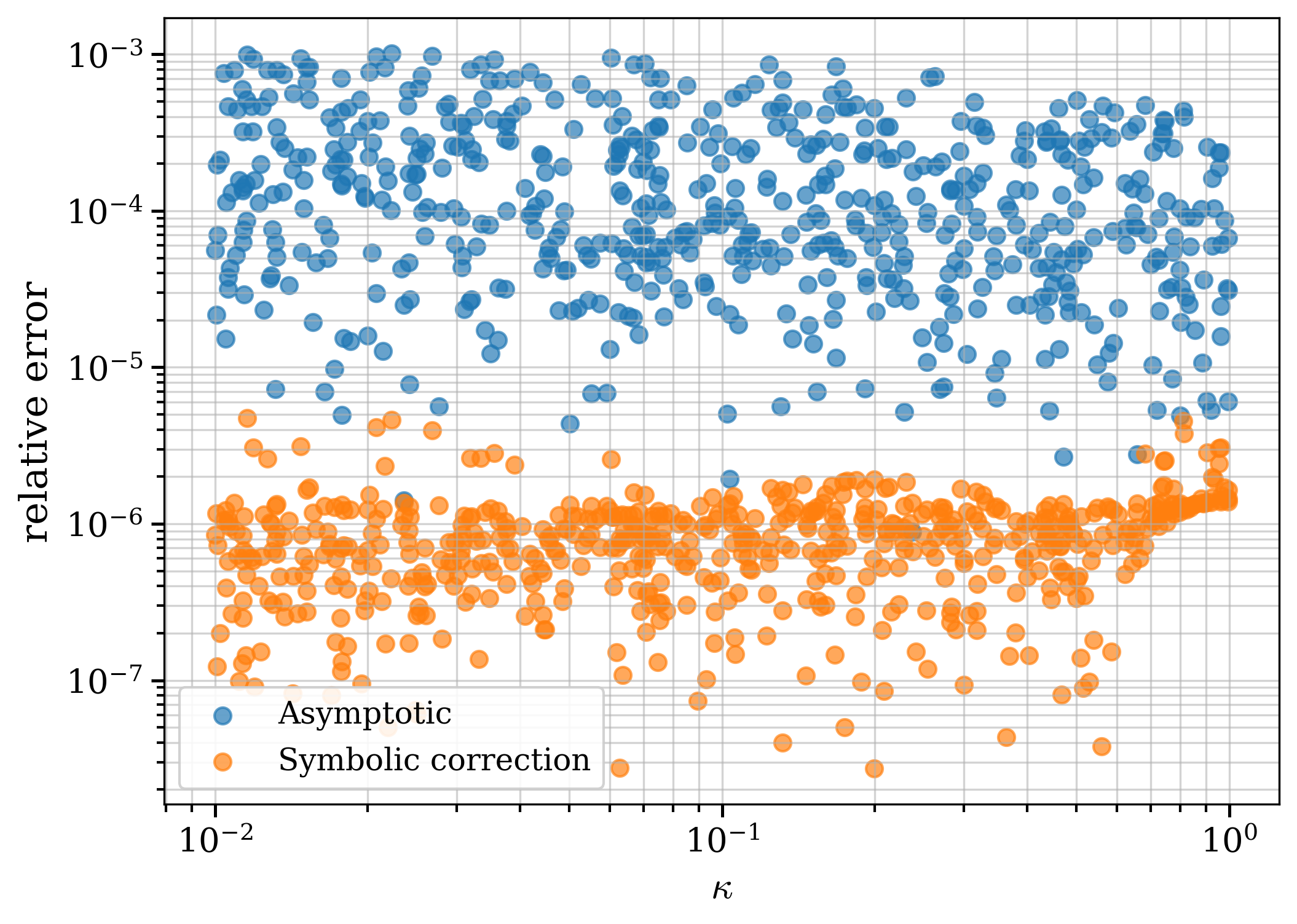}
        \caption{Second dimer branch: error of the compact symbolic-regression correction.}
        \label{fig:dimer_symbolic_branch2}
    \end{subfigure}
    \caption{Symbolic-regression correction for the two dimer branches as functions of the separation distance $\kappa$. The symbolic formulas preserve the accuracy of the learned residual correction while providing explicit expressions involving the asymptotic eigenvalues, branch splitting, and the expected powers of $\delta$ and $q_{j}$.}
    \label{fig:dimer:symbolic correction}
\end{figure}

We then evaluate the symbolic-regression formulas for the two dimer branches. Figures~\ref{fig:dimer_symbolic_branch1} and~\ref{fig:dimer_symbolic_branch2} show that the symbolic corrections preserve the accuracy of the learned residual model while providing explicit expressions depending on the asymptotic eigenvalues, the branch splitting, and the expected powers of $\delta$ and $q_{j}$.

The improvement factors are smaller than in the single-resonator case. This is expected: the dimer residual contains several sources of error at once, including the self-interaction error inherited from the single-resonator approximation, the error in the interaction blocks of the reduced operator, and the sensitivity of the two resonance branches to the eigenvalue splitting. Moreover, small changes in $\kappa$ can modify the coupling strength substantially. The dimer results should therefore be interpreted as a more demanding interaction test rather than as a direct confirmation of the single-resonator behaviour.

We also examine the dimer errors as functions of the size parameter $\delta$. While the separation distance $\kappa$ controls the strength of the interaction between the two resonators, the parameter $\delta$ is the small parameter in the subwavelength asymptotic regime. This provides a complementary view to the previous plots: the $\kappa$-dependence describes the coupling effects, whereas the $\delta$-dependence shows how the approximation behaves with respect to the asymptotic limit.

\begin{figure}[ht]
    \centering
    \begin{subfigure}[t]{0.47\textwidth}
        \centering
        \includegraphics[width=\linewidth]{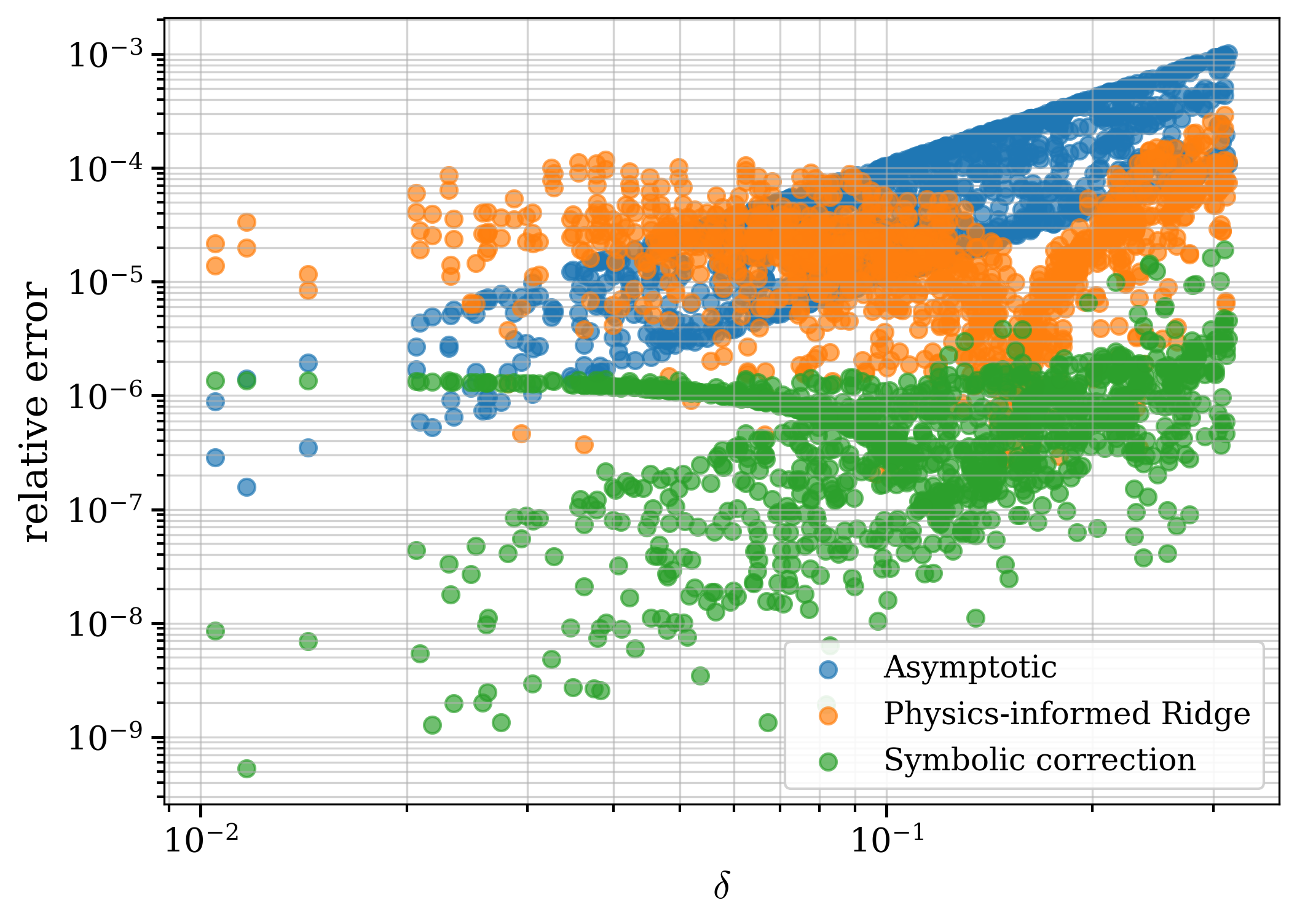}
        \caption{Relative error versus the size parameter $\delta$.}
        \label{fig:dimer_error_delta}
    \end{subfigure}
    \hfill
    \begin{subfigure}[t]{0.43\textwidth}
        \centering
        \includegraphics[width=\linewidth]{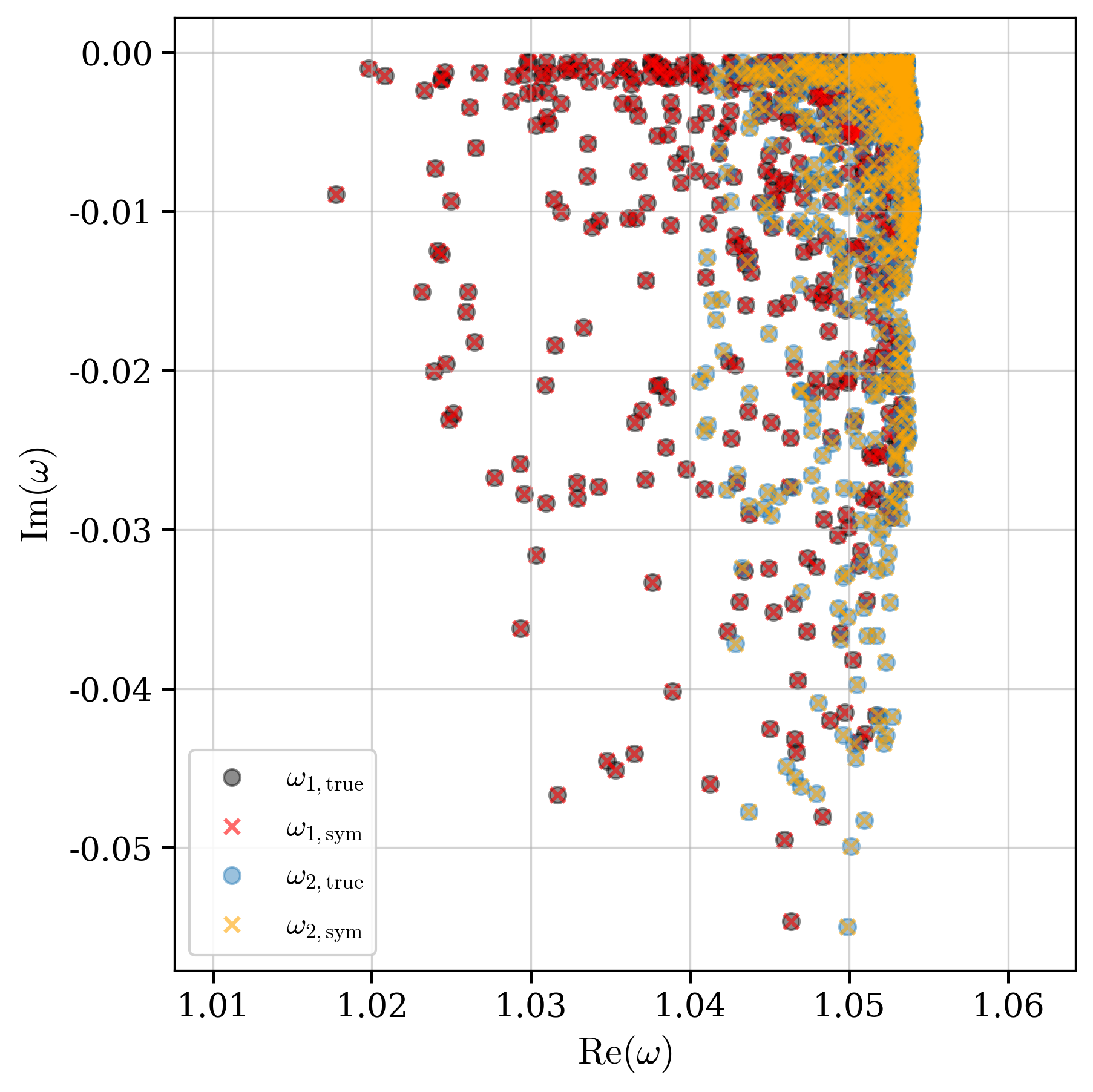}
        \caption{Reconstruction of both dimer branches in the complex plane.}
        \label{fig:dimer_symbolic_complex}
    \end{subfigure}
    \caption{Dimer prediction results. Left: relative error of the asymptotic approximation and of the corrected models as a function of the asymptotic parameter $\delta$, with both resonance branches included. Right: symbolic-regression reconstruction of the two resonance branches in the complex plane. The corrected models reduce the baseline error across the sampled range, while the broader scatter compared with the single-resonator case reflects the additional dependence on the separation distance $\kappa$, damping and branch coupling.}
    \label{fig:dimer:error_vs_delta}
\end{figure}

Figure~\ref{fig:dimer_error_delta} shows the relative errors for both dimer branches as functions of $\delta$. The baseline asymptotic approximation exhibits a clear increase in error as $\delta$ grows, consistently with the fact that the asymptotic expansion is most accurate in the small-$\delta$ regime. The physics-informed Ridge correction reduces the error across the full range of $\delta$, although the corrected errors display a broader scatter than in the single-resonator case. This is expected, since in the dimer configuration the residual depends not only on $\delta$, but also on the separation distance $\kappa$, the damping parameter $\gamma$ and the resonance branch.

The symbolic correction remains systematically below the Ridge correction over the sampled range of $\delta$, showing that the symbolic formulas capture additional structure in the residual. The remaining scatter reflects the fact that the dimer residual is genuinely multi-parameter: samples with similar values of $\delta$ may correspond to different coupling strengths and therefore to different correction magnitudes. In Figure~\ref{fig:dimer_symbolic_complex}, the joint reconstruction of the two branches in the complex plane is shown. The symbolic predictions align closely with the reference resonances for both branches. Quantitatively, the mean relative error is reduced from approximately $1.22\times 10^{-4}$ for the asymptotic model to approximately $8.35\times 10^{-7}$ for the compact symbolic formulas.

\subsection{Ablation study}

We now assess the role of feature design more systematically. The purpose of this ablation study is to determine whether the improvement comes merely from increasing the number of input variables, or from adding variables that reflect the asymptotic structure of the resonance problem.

We compare the following nested feature sets:
\begin{itemize}
    \item $M0$: raw parameters only, namely $(\delta,\gamma)$ for the single resonator and $(\delta,\gamma,\kappa)$ for the dimer;
    \item $M1$: raw parameters together with logarithmic transforms;
    \item $M2$: asymptotic resonance or spectral quantities, such as $\omega_{\mathrm{asymp}}$, $\lambda_{\mathrm{asymp}}$, and, for the dimer, the branch splitting;
    \item $M3$: algebraic asymptotic scales, such as $\delta^2$, $\delta^2\log\delta$, and $\gamma\delta^2$;
    \item $M4$: frequency-dependent logarithmic scales, such as $\delta^2 q^2$, $\delta^2 q^2\log\delta$, and $\delta^2 q^2\log q$.
\end{itemize}

The results are reported in Table~\ref{tab:enhanced_ablation}. For the single resonator, raw parameters already improve the baseline, reducing the mean relative error from $9.09\times10^{-4}$ to $1.28\times10^{-4}$. Adding logarithmic transforms and the asymptotic resonance further improves the prediction. However, the dominant gain occurs when the algebraic asymptotic scales are introduced: the mean error drops from $8.23\times10^{-6}$ for $M2$ to $1.82\times10^{-7}$ for $M3$. The full feature set $M4$, which includes the frequency-dependent logarithmic scales, gives the lowest mean error, $1.25\times10^{-7}$.

The same hierarchy is observed for the dimer. The raw and logarithmic parameter sets $M0$ and $M1$ provide moderate improvements, while the inclusion of asymptotic spectral quantities in $M2$ further reduces the error. The main reduction again occurs when algebraic asymptotic scales are added: the mean error decreases from $2.21\times10^{-5}$ for $M2$ to $3.67\times10^{-6}$ for $M3$. Adding the frequency-dependent logarithmic scales in $M4$ gives an additional improvement, reaching $2.18\times10^{-6}$.

This progression shows that the correction is not improved merely because more variables are supplied to the regression model. Rather, the improvement is driven by asymptotically meaningful information. In particular, the sharp decrease from $M2$ to $M3$, followed by the additional improvement from $M3$ to $M4$, supports the claim that the learned residual is organized by the same scales that appear in the asymptotic expansion.

\begin{table}[ht]
    \centering
    \begin{tabular}{llccccc}
        \hline
        Configuration & Model & Features & Median & 90th perc. \\
        \hline
        Single & Asymptotic baseline & 0 & $5.85\times10^{-4}$ & $2.28\times10^{-3}$ \\
        Single & $M0$: raw & 2 & $1.26\times10^{-4}$ & $2.12\times10^{-4}$ \\
        Single & $M1$: raw + logs & 4 & $3.81\times10^{-5}$ & $6.38\times10^{-5}$ \\
        Single & $M2$: + asymptotic resonance & 6 & $6.95\times10^{-6}$ & $1.49\times10^{-5}$ \\
        Single & $M3$: + algebraic scales & 9 & $1.42\times10^{-7}$ & $3.73\times10^{-7}$ \\
        Single & $M4$: + frequency-log scales & 12 & $8.19\times10^{-8}$ & $2.61\times10^{-7}$ \\
        \hline
        Dimer & Asymptotic baseline & 0 & $5.16\times10^{-5}$ & $3.43\times10^{-4}$ \\
        Dimer & $M0$: raw & 3 & $2.77\times10^{-5}$ & $8.19\times10^{-5}$ \\
        Dimer & $M1$: raw + logs & 6 & $1.94\times10^{-5}$ & $7.19\times10^{-5}$ \\
        Dimer & $M2$: + asymptotic spectrum & 13 & $1.23\times10^{-5}$ & $5.40\times10^{-5}$ \\
        Dimer & $M3$: + algebraic scales & 16 & $2.59\times10^{-6}$ & $7.74\times10^{-6}$ \\
        Dimer & $M4$: + frequency-log scales & 22 & $1.54\times10^{-6}$ & $4.68\times10^{-6}$ \\
        \hline
    \end{tabular}
    \caption{Ablation study for the Ridge residual correction. The table reports the distribution of relative errors on the held-out test sets. The nested feature hierarchy shows that the main improvement is obtained when asymptotic and physics-informed scales are added, rather than merely by increasing the number of input variables.}
    \label{tab:enhanced_ablation}
\end{table}

\subsection{Stability of the symbolic expressions}

Symbolic regression involves a stochastic search over a large space of candidate expressions. Therefore, a single symbolic formula should not be overinterpreted as a unique correction law. To assess the stability of the symbolic structure, we repeat the assisted symbolic-regression procedure over ten independent runs with different random seeds. For each run, symbolic regression is applied separately to the real and imaginary parts of the residual. In the dimer case, it is also applied separately to the two resonance branches.

We focus on the two assisted strategies introduced above. In the first one, called \emph{asymptotic preconditioning}, a next-order proxy is first fitted using the asymptotic feature family $\delta^2 q_{\mathrm{asymp}}^2$, $\delta^2 q_{\mathrm{asymp}}^2\log\delta$, $\delta^2 q_{\mathrm{asymp}}^2\log q_{\mathrm{asymp}}$ and symbolic regression is then applied to the remaining residual. In the second one, called \emph{asymptotic feature augmentation}, the proxy is supplied as an additional input variable and symbolic regression is free to use it or ignore it.

For each selected expression, we record whether the proxy, the asymptotic logarithmic feature family, and the main remaining symbolic structures appear in the formula. The results are summarized in Table~\ref{tab:assisted_symbolic_stability}. For the single resonator, the counts for symbolic expressions are taken over $20$ formulas, corresponding to $10$ runs times real and imaginary parts. For the dimer, they are taken over $40$ formulas, corresponding to $10$ runs times two branches times real and imaginary parts.

\begin{table}[ht]
    \centering
    \begin{tabular}{llcc}
        \hline
        Strategy & Feature or structure & Single & Dimer \\
        \hline
        \multirow{4}{*}{Feature augmentation}
        & proxy & $14/20$ & $15/40$ \\
        & $\delta^2 q^2$ & $14/20$ & $21/40$ \\
        & $\delta^2 q^2\log\delta$ & $1/20$ & $3/40$ \\
        & $\delta^2 q^2\log q$ & $6/20$ & $13/40$ \\
        \hline
        \multirow{4}{*}{Preconditioning}
        & proxy & $20/20$ & $40/40$ \\
        & $\delta^2 q^2$ & $20/20$ & $40/40$ \\
        & $\delta^2 q^2\log\delta$ & $20/20$ & $40/40$ \\
        & $\delta^2 q^2\log q$ & $20/20$ & $40/40$ \\
        \hline
        \multirow{4}{*}{Remainder structure}
        & $\delta^2$-type terms & $20/20$ & $40/40$ \\
        & $\gamma\delta^2$-type terms & $8/20$ & $8/40$ \\
        & $\omega_{\mathrm{asymp}}$-dependent terms & $12/20$ & -- \\
        & spectral/splitting terms & -- & $15/40$, $9/40$ \\
        \hline
    \end{tabular}
    \caption{Stability of the assisted symbolic-regression expressions. The counts are computed over ten independent runs, corresponding to \(20\) symbolic expressions for the single resonator and \(40\) expressions for the dimer. The table reports how often each feature family or symbolic structure appears in the selected expressions. In the dimer case, the spectral/splitting row reports the frequencies of \(\lambda_{\mathrm{asymp}}\)-dependent terms and splitting-dependent terms, respectively.}
    \label{tab:assisted_symbolic_stability}
\end{table}

The table distinguishes two notions of stability. In the preconditioned setting, the next-order proxy is present by construction. Therefore, the logarithmic feature family $\delta^2 q_{\mathrm{asymp}}^2\log\delta$, $\delta^2 q_{\mathrm{asymp}}^2\log q_{\mathrm{asymp}}$ appears in every run. This confirms that the preconditioned symbolic correction consistently incorporates the asymptotically motivated next-order scale. The remaining symbolic correction is then fitted on top of this proxy, and its most stable structure is the $\delta^2$-type contribution, which appears in every run for both configurations.

The feature-augmentation setting is more selective, because PySR is given the proxy as an input but is not forced to use it. In this case, the proxy is selected in $14/20$ single-resonator expressions and $15/40$ dimer expressions. The algebraic frequency-dependent scale $\delta^2 q^2$ is also selected frequently. The individual logarithmic variables are selected less often, especially $\delta^2 q^2\log\delta$. This indicates that, over the sampled parameter range, symbolic regression may represent part of the logarithmic contribution through the proxy or through simpler algebraic combinations. Thus the logarithmic terms are stable in the asymptotically preconditioned correction, while unconstrained selection within the augmented feature set tends to favor more compact representatives.

The prediction errors across independent runs are reported in Table~\ref{tab:assisted_symbolic_stability_errors}. The feature-augmentation strategy gives the most stable accuracy, with mean relative errors of order $10^{-7}$ for the single resonator and $10^{-6}$ for the dimer. The preconditioned strategy has a much smaller median error than mean error, indicating that most runs perform well but a few runs produce less accurate symbolic remainders. The proxy alone is less accurate, which shows that the asymptotic next-order family captures an important part of the residual but still benefits from a learned symbolic correction.

\begin{table}[ht]
    \centering
    \begin{tabular}{llccc}
        \hline
        Configuration & Strategy & Mean & Median & Standard deviation \\
        \hline
        Single & next-order proxy & $7.31\times10^{-3}$ & $1.85\times10^{-5}$ & $1.97\times10^{-2}$ \\
        Single & asymptotic preconditioning & $5.07\times10^{-5}$ & $1.54\times10^{-6}$ & $1.22\times10^{-4}$ \\
        Single & feature augmentation & $1.66\times10^{-7}$ & $7.53\times10^{-8}$ & $1.53\times10^{-7}$ \\
        Dimer & next-order proxy & $3.45\times10^{-3}$ & $9.68\times10^{-5}$ & $5.26\times10^{-3}$ \\
        Dimer & asymptotic preconditioning & $1.10\times10^{-4}$ & $1.86\times10^{-5}$ & $1.76\times10^{-4}$ \\
        Dimer & feature augmentation & $4.87\times10^{-6}$ & $3.18\times10^{-6}$ & $4.11\times10^{-6}$ \\
        \hline
    \end{tabular}
    \caption{Prediction errors for the assisted symbolic-regression stability study over ten independent runs. The reported quantities are computed from the mean relative error of each run.}
    \label{tab:assisted_symbolic_stability_errors}
\end{table}

Overall, the stability study supports the interpretation used throughout the paper. The asymptotic logarithmic family is not merely an arbitrary collection of input variables: it provides a stable next-order proxy when the symbolic regression is guided by the asymptotic structure. At the same time, the individual symbolic formulas are not unique, and PySR may choose different but comparably accurate algebraic representations across runs. This reinforces the view that symbolic regression should be interpreted here as a tool for discovering compact data-driven correction laws consistent with the asymptotic structure, rather than as a rigorous derivation of the next asymptotic term.

\subsection{Assisted symbolic regression}

Finally, we test the assisted symbolic-regression strategy described in Section~\ref{sec:symbolic_regression}. The goal is to help symbolic regression by explicitly providing information about the expected next-order asymptotic structure. We first construct an asymptotic proxy for the next-order correction using the feature family $\delta^2 q_{\mathrm{asymp}}^2\log\delta$, $\delta^2 q_{\mathrm{asymp}}^2\log q_{\mathrm{asymp}}$, $\delta^2 q_{\mathrm{asymp}}^2$. We then compare three strategies: using this proxy alone, asking symbolic regression to learn only the remaining discrepancy after subtracting the proxy, and giving the proxy as an additional input variable to symbolic regression. These correspond respectively to the next-order proxy, asymptotic preconditioning, and asymptotic feature augmentation.

\begin{figure}[ht]
    \centering
    \begin{subfigure}[t]{0.48\textwidth}
        \centering
        \includegraphics[width=\linewidth]{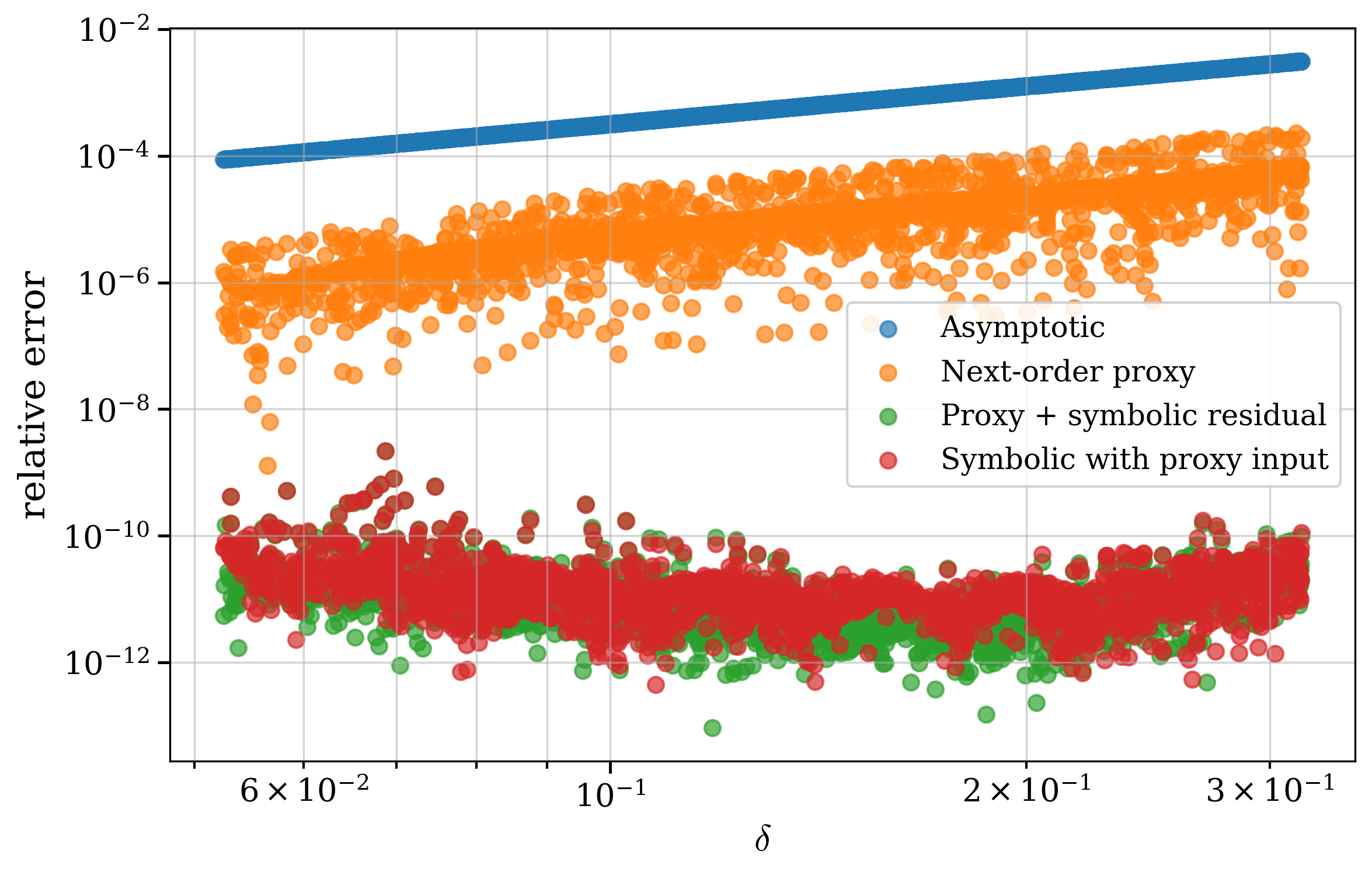}
        \caption{Single resonator.}
        \label{fig:single_advisor_symbolic}
    \end{subfigure}
    \hfill
    \begin{subfigure}[t]{0.48\textwidth}
        \centering
        \includegraphics[width=\linewidth]{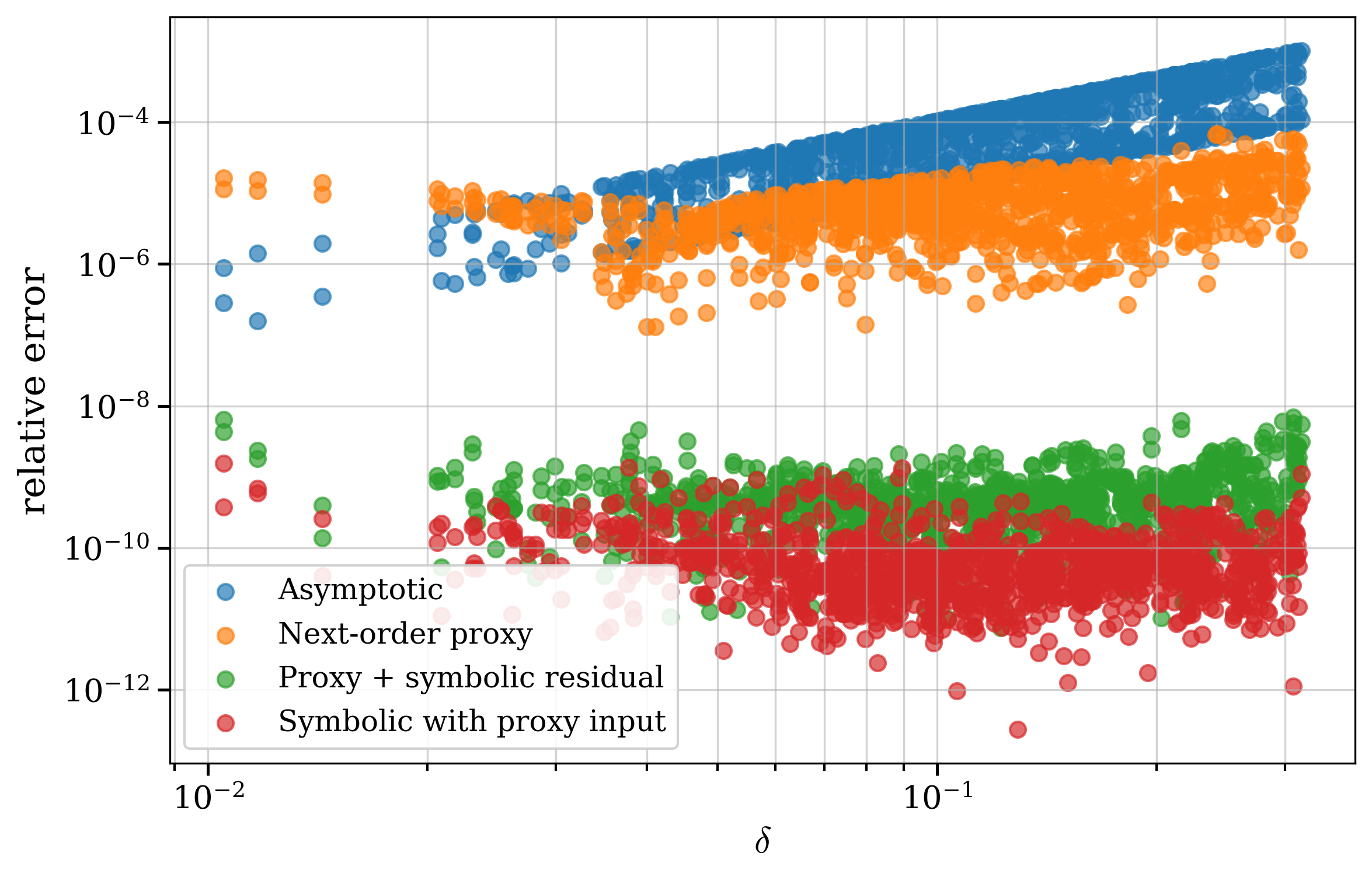}
        \caption{Dimer, with both resonance branches included.}
        \label{fig:dimer_advisor_symbolic_delta}
    \end{subfigure}
    \caption{Assisted symbolic regression using a next-order asymptotic proxy. The proxy already improves the asymptotic or reduced approximation, and symbolic regression further reduces the remaining residual by several orders of magnitude. The dimer plot is more dispersed because the correction also depends on the separation distance and on branch splitting.}
    \label{fig:advisor_symbolic_delta}
\end{figure}

Figure~\ref{fig:advisor_symbolic_delta} shows the resulting errors as functions of the size parameter $\delta$, for the single resonator and for the dimer. In the single-resonator case, the asymptotic proxy already improves the baseline approximation, reducing the mean relative error from approximately $9.09\times10^{-4}$ to approximately $1.83\times10^{-5}$. Symbolic regression on the remaining correction further reduces the error to approximately $1.89\times10^{-11}$, while giving the proxy as an additional input gives approximately $2.22\times10^{-11}$. Thus, once the next-order asymptotic structure is made explicit, symbolic regression is able to recover almost all of the remaining residual.

The same behaviour is observed for the dimer, although the error cloud is broader. This is expected because the dimer residual depends not only on $\delta$, but also on the separation distance $\kappa$, on the damping parameter, and on the splitting between the two branches. In this case, the asymptotic proxy reduces the mean relative error from approximately $1.22\times10^{-4}$ to approximately $8.27\times10^{-6}$. Symbolic regression on the remaining correction gives approximately $6.04\times10^{-10}$, while using the proxy as an input feature gives approximately $9.30\times10^{-11}$. The feature-augmentation strategy therefore gives the lowest mean error in the dimer case.

These results show that symbolic regression is most effective when it is guided by the asymptotic structure. The next-order proxy captures a substantial part of the missing correction, while the symbolic step learns the remaining discrepancy. The comparison between asymptotic preconditioning and asymptotic feature augmentation also shows that the asymptotic information can be injected either through the target or through the input representation. Both strategies perform very well for the single resonator, while feature augmentation is slightly more effective for the dimer. This supports the interpretation that the learned residual is governed by the same scales as the neglected asymptotic terms, even in the presence of coupling and branch splitting.

\subsection{Robustness with respect to noisy training data}

We now investigate the robustness of the correction models with respect to noise in the reference resonances. This experiment is motivated by the fact that, in practical settings, the values of $\omega_{\mathrm{ref}}$ may come from numerical solvers or measurements and may therefore contain errors.

\begin{figure}[ht]
    \centering
    \begin{subfigure}[t]{0.47\textwidth}
        \centering
        \includegraphics[width=\linewidth]{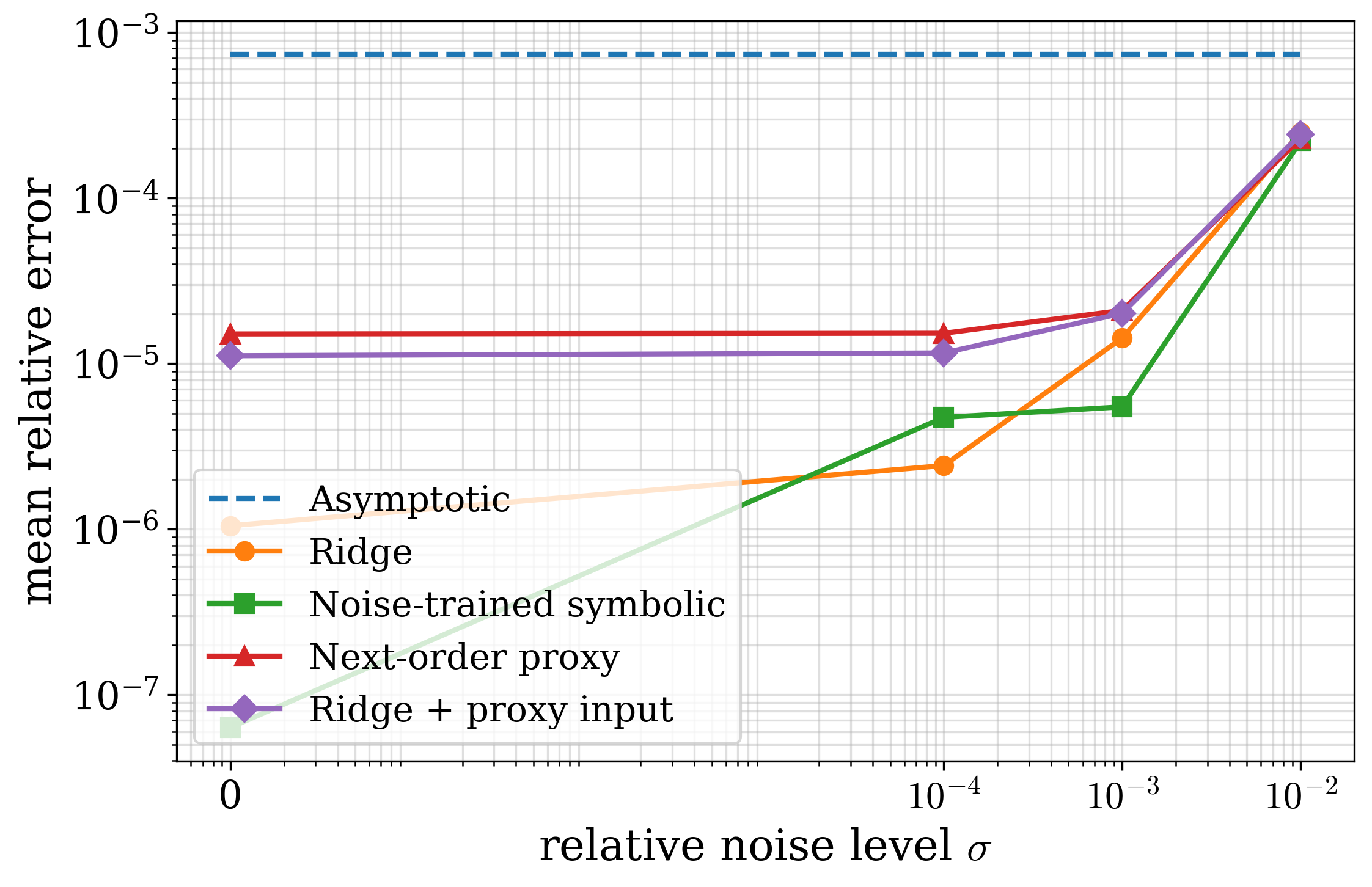}
        \caption{Single resonator.}
        \label{fig:single_noise_robustness}
    \end{subfigure}
    \hfill
    \begin{subfigure}[t]{0.47\textwidth}
        \centering
        \includegraphics[width=\linewidth]{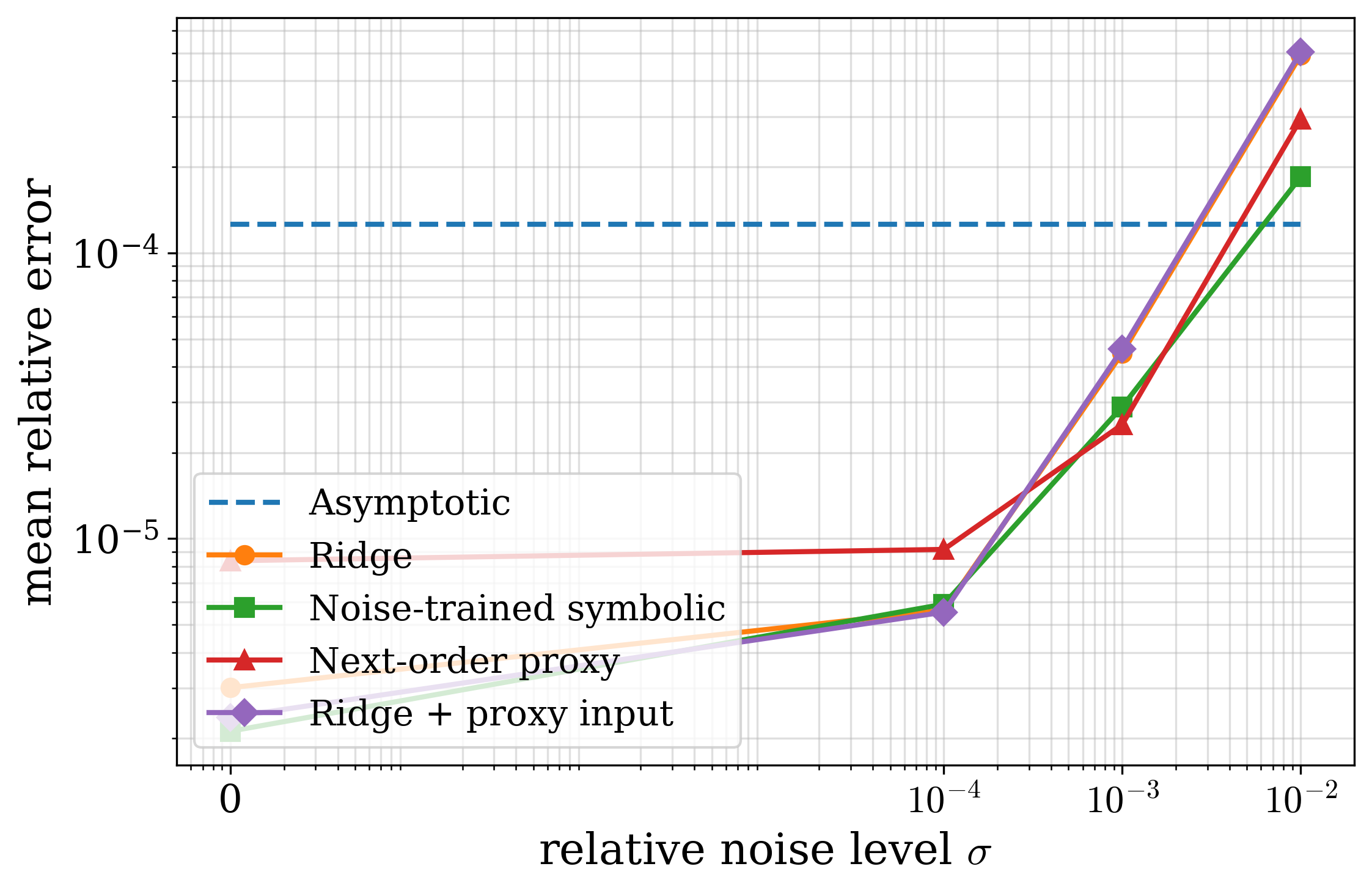}
        \caption{Dimer.}
        \label{fig:dimer_noise_robustness}
    \end{subfigure}
    \caption{Robustness with respect to noisy training resonances. The corrections remain stable for moderate noise levels and deteriorate when the noise becomes comparable to the residual.}
    \label{fig:noise_robustness}
\end{figure}

We also test the effect of perturbing the reference resonances used for training. Figure~\ref{fig:noise_robustness} shows that the corrections remain stable for moderate noise levels and degrade once the noise becomes comparable to the residual being learned. The dimer is more sensitive, as expected from its branch-dependent coupling corrections. The precise noise model, noise levels, and dataset sizes are given in Appendix~\ref{app:numerical_details}.

A clear transition occurs between $\sigma=10^{-3}$ and $\sigma=10^{-2}$. For $\sigma \leq 10^{-3}$, the corrected models remain well below the asymptotic error, showing that the residual-learning strategy is robust to moderate perturbations of the training data. At $\sigma=10^{-2}$, the noise dominates the residual and the learned corrections deteriorate substantially. This behaviour is expected: once the noise level is larger than the correction itself, the training target no longer contains enough reliable information to identify the residual accurately.

The single-resonator case is more stable than the dimer case. In the dimer, the residual contains both self-interaction corrections and branch-dependent coupling effects, making the learning problem more sensitive to perturbations. Nevertheless, the same qualitative behaviour is observed in both configurations. The noise-trained symbolic correction follows the same degradation pattern as the Ridge and proxy-based corrections, confirming that symbolic regression can also be retrained from noisy residuals when sufficiently large datasets are available. The next-order proxy and the Ridge model augmented with the proxy feature behave similarly, reflecting the fact that both are constrained by the same asymptotic structure.

\subsection{Effect of the number of training measurements}

We finally investigate how the robustness of the different correction strategies depends on the number of available training measurements. For each training-set size, the models are trained under noisy resonances and evaluated on a fixed held-out test set. The same relative noise levels as before are considered, namely $\sigma \in \{0,10^{-4},10^{-3},10^{-2}\}$. This experiment complements the previous noise study by separating the effect of the noise amplitude from the effect of the amount of training data.

\begin{figure}[ht]
    \centering
    \includegraphics[width=0.95\textwidth]{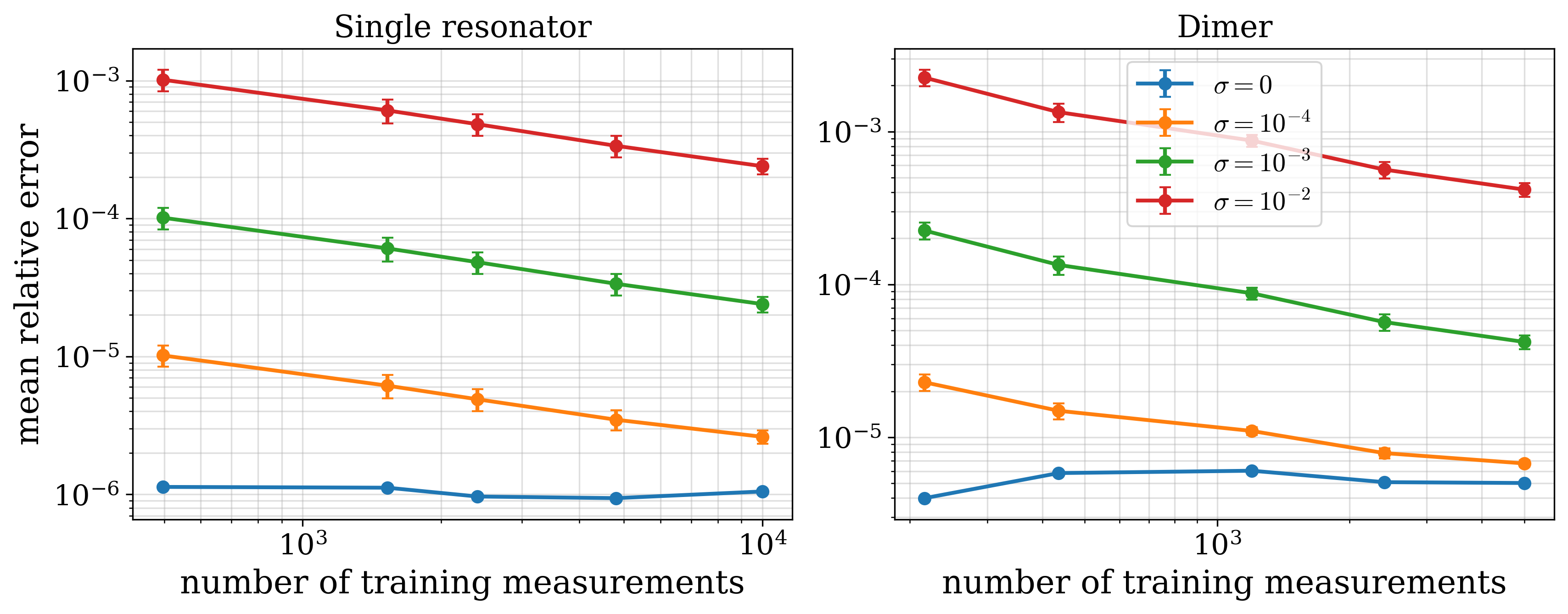}
    \caption{Effect of the number of training measurements on the Ridge residual correction. The error decreases as the amount of training data increases, especially for noisy measurements.}
    \label{fig:training_size_ridge}
\end{figure}

Figure~\ref{fig:training_size_ridge} shows the behaviour of the Ridge residual correction. Increasing the number of training measurements significantly improves the prediction under noisy data. The improvement is especially visible for the larger noise levels, where averaging over more samples stabilizes the learned residual. For clean data, the error is already close to its limiting value and therefore changes only mildly with the training size.

\begin{figure}[ht]
    \centering
    \includegraphics[width=0.95\textwidth]{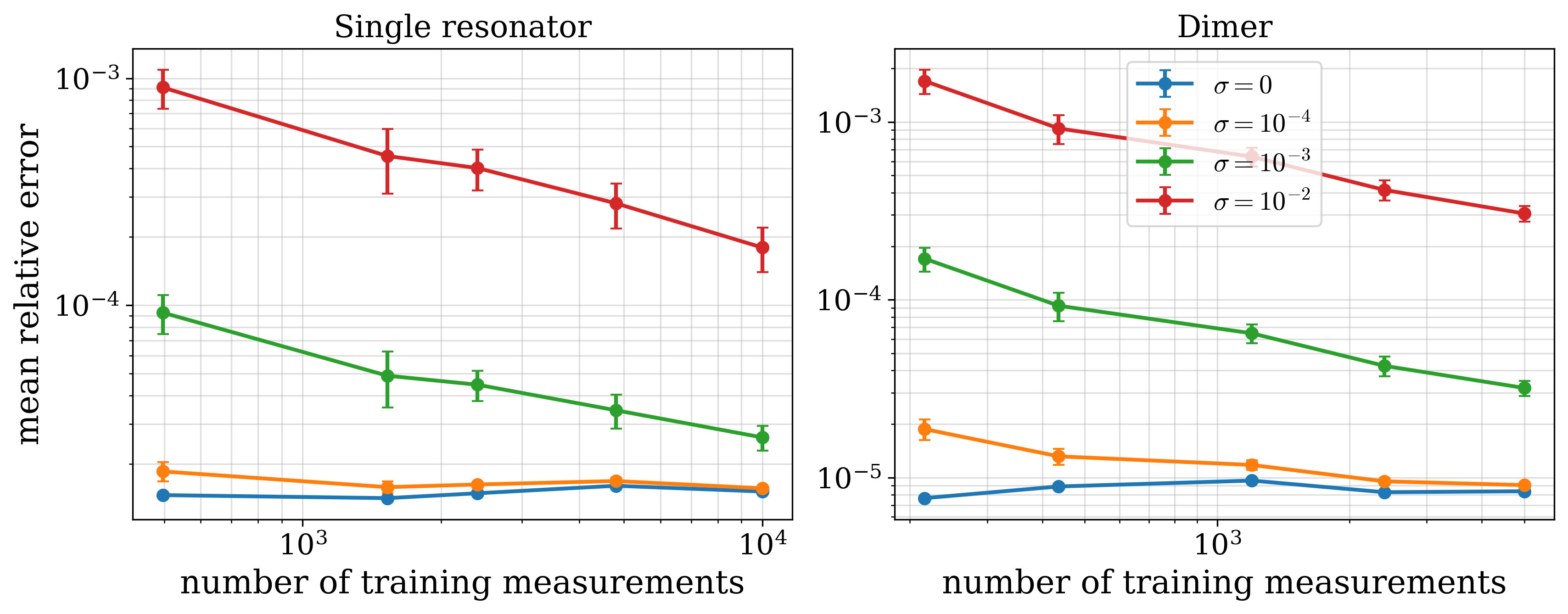}
    \caption{Effect of the number of training measurements on the next-order asymptotic proxy. Increasing the training size improves the stability of the calibrated proxy under noisy measurements.}
    \label{fig:training_size_proxy}
\end{figure}

The corresponding result for the next-order proxy is shown in Figure~\ref{fig:training_size_proxy}. The same qualitative behaviour is observed: the proxy becomes more stable as the number of measurements increases. This is consistent with the fact that the proxy is calibrated from noisy residual data using asymptotically motivated features. The effect is again most pronounced for the larger noise levels.

\begin{figure}[ht]
    \centering
    \includegraphics[width=0.95\textwidth]{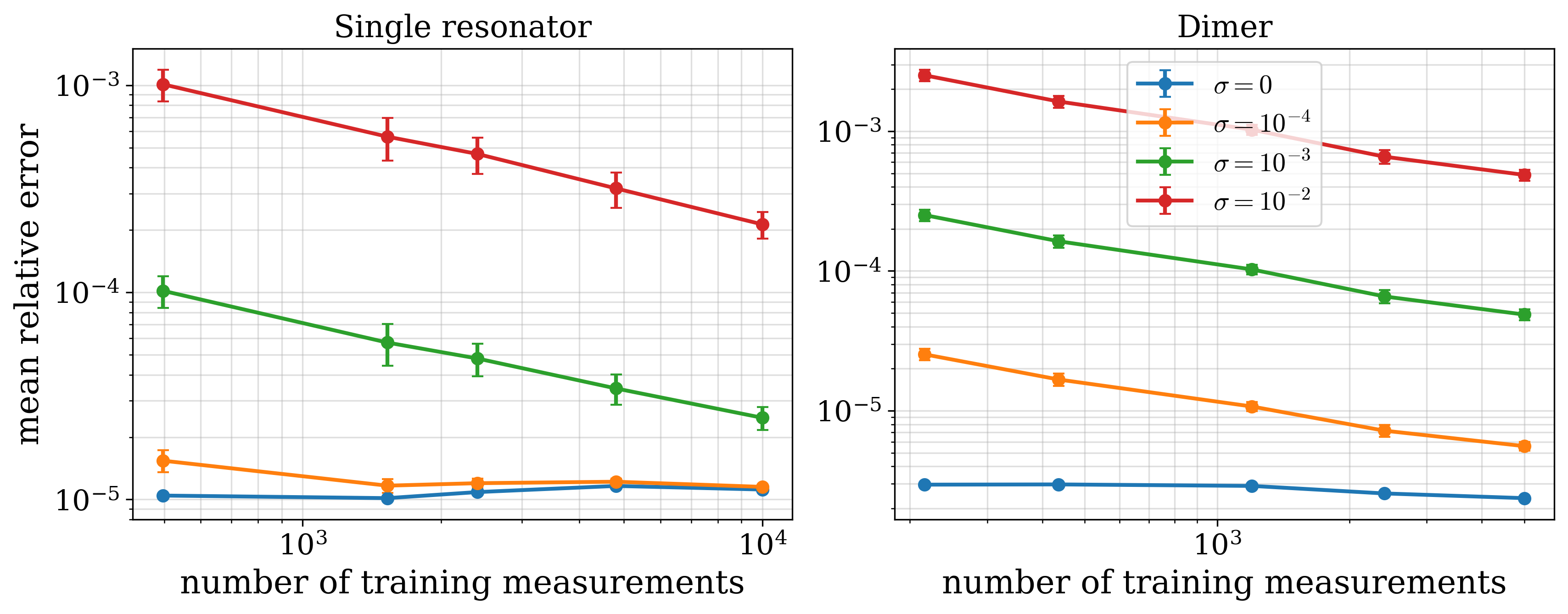}
    \caption{Effect of the number of training measurements on the Ridge correction assisted by the next-order proxy. The proxy feature improves the stability of the learned correction, with a particularly visible gain in the dimer case.}
    \label{fig:training_size_ridge_proxy}
\end{figure}

Figure~\ref{fig:training_size_ridge_proxy} reports the result obtained when the next-order proxy is supplied as an additional input to the Ridge model. This assisted correction is also improved by larger training sets. In the dimer case, the gain is particularly clear, reflecting the fact that the coupled problem requires more data to stabilize the branch-dependent correction.

\begin{figure}[ht]
    \centering
    \includegraphics[width=0.95\textwidth]{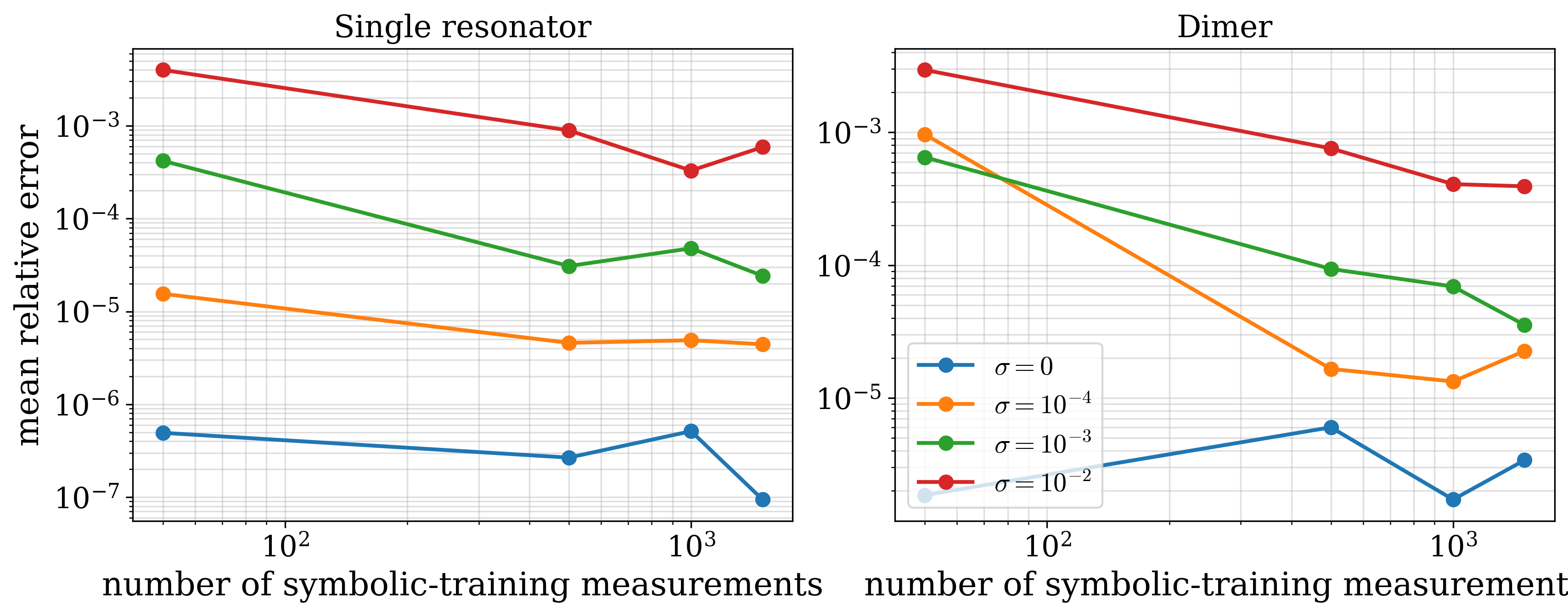}
    \caption{Effect of the number of symbolic-training measurements on PySR rediscovery under noisy data. For each training size and noise level, symbolic regression is rerun from the noisy residuals. Larger training sets improve the robustness of the rediscovered symbolic correction, especially at higher noise levels.}
    \label{fig:training_size_symbolic}
\end{figure}

Finally, Figure~\ref{fig:training_size_symbolic} shows the same experiment for symbolic regression, where PySR is retrained independently for each training size and each noise level. Increasing the symbolic-training set improves the rediscovered formulas, especially for $\sigma=10^{-3}$ and $\sigma=10^{-2}$. The curves are not expected to be perfectly monotone, since symbolic regression solves a non-convex model-selection problem and each run may discover a different analytic expression. Nevertheless, the overall trend shows that larger datasets make symbolic rediscovery more stable under noise.

Overall, these experiments confirm that the effect of measurement noise can be mitigated by increasing the number of training samples. Ridge-based corrections benefit from this most directly, while symbolic regression also becomes more stable when more data are available, although its non-convex search may lead to mild non-monotonicity.

\subsection{Summary}

Across both the single-resonator and dimer configurations, the proposed corrections substantially improve the asymptotic or reduced approximations. Ridge regression provides a stable residual model, while symbolic regression yields compact analytic formulas that retain high accuracy and expose the structure of the correction. The ablation studies, now performed for both configurations, show that this improvement is not simply due to adding more variables, but relies on features that reflect the asymptotic structure of the problem. In particular, the logarithmic frequency-dependent features $\delta^2 q_{\mathrm{asymp}}^2\log\delta$ and $\delta^2 q_{\mathrm{asymp}}^2\log q_{\mathrm{asymp}}$ play an important role in capturing the next-order behaviour.

The assisted symbolic-regression experiments further confirm this interpretation. For both the single resonator and the dimer, an explicit next-order proxy already reduces the asymptotic error, and symbolic regression then captures the remaining discrepancy. In the dimer case, the error distribution is broader because the correction depends not only on $\delta$, but also on the separation distance $\kappa$ and on the splitting between the two resonant branches. Nevertheless, the same hierarchy of models is observed, showing that the learned residual remains organized by the expected asymptotic scales.

The noise experiments show that the residual-learning framework is robust to moderate perturbations of the training resonances, especially when the number of training measurements is increased. Ridge and proxy-based corrections degrade gradually as the noise level grows, and remain effective as long as the noise is not larger than the residual being learned. The dimer case is more sensitive than the single-resonator case, reflecting the additional difficulty of learning coupled, branch-dependent corrections.

Finally, the noisy-test experiment distinguishes between using a symbolic formula discovered from clean data and rediscovering a symbolic formula from noisy data. In both the single-resonator and dimer cases, the clean symbolic formula performs slightly better. This suggests that PySR is most effective as a tool for discovering clean structural corrections, rather than as a procedure for fitting noisy measurements directly.

\section{Discussion and perspectives}\label{sec:Conclusion and perspectives}

In this work, we developed an asymptotics-guided learning framework for resonance prediction in dispersive resonator systems. The central point is that asymptotic analysis is useful not only because it provides a leading-order approximation, but also because it identifies the variables and scales that organize the error of this approximation. We therefore use the asymptotic model as a baseline and learn only the residual between this prediction and the numerical reference resonance.

The numerical results show that this residual is highly structured. For both the single resonator and the dimer, physics-informed features derived from the subwavelength expansion substantially improve the prediction compared with raw-parameter learning. In particular, the logarithmic frequency-dependent scales suggested by the two-dimensional operator expansion provide effective features for correcting the missing higher-order terms. The dimer case is more demanding, since the residual also contains branch-dependent interaction effects, but the same asymptotics-guided strategy remains effective.

Symbolic regression adds an interpretable layer to this correction procedure. The formulas obtained should not be interpreted as rigorous asymptotic expansions; they are data-driven correction laws valid in the sampled parameter regime. Their significance is that they recover compact expressions organized around the asymptotic scales, including size-dependent, frequency-dependent, damping-dependent, and branch-splitting terms. The assisted symbolic-regression and stability experiments further support the view that symbolic regression is most effective when guided by the asymptotic structure rather than used as a purely black-box search.

Overall, the results suggest a general methodology: use asymptotic analysis to identify the structure of the residual, use learning to correct the missing terms, and use symbolic regression to extract compact interpretable formulas. This provides a bridge between asymptotic modelling and data-driven prediction for nonlinear spectral problems in dispersive media.

\appendix

\section{Numerical and learning details}\label{app:numerical_details}

This appendix summarizes the numerical choices used to generate the datasets and reproduce the learning experiments. All computations are implemented in Python. Data handling uses \texttt{NumPy} and \texttt{pandas}, Ridge regression and feature standardization use \texttt{scikit-learn}, and symbolic regression is performed with PySR, interfaced with \texttt{SymbolicRegression.jl}. Random seeds are fixed throughout the experiments.

\subsection{Parameter ranges and dataset construction}

All parameter values in this appendix are non-dimensional. For the single-resonator experiments, the size parameter and damping parameter are sampled logarithmically $\delta \sim 10^{U(-2,-1/2)}$ and $\gamma \sim 10^{U(-3,-1)}$,
where $U(a,b)$ denotes the uniform distribution on the interval $[a,b]$. Thus, $\delta \in [10^{-2},10^{-1/2}]$ and $\gamma \in [10^{-3},10^{-1}]$.
In the final single-resonator dataset used for the main plots, the realized range is approximately $\delta \in [5.25\times 10^{-2},3.16\times 10^{-1}]$ and $\gamma \in [1.00\times 10^{-3},9.998\times 10^{-2}]$.
The lower end of the realized $\delta$-range is larger than $10^{-2}$ because samples outside the retained physical branch are discarded.

For the dimer experiments, the size parameter, damping parameter, and separation distance are sampled as $\delta \sim 10^{U(-2,-1/2)}$, $\gamma \sim 10^{U(-3,-1)}$ and $\kappa \sim 10^{U(-2,0)}$.
Hence, $\kappa \in [10^{-2},1]$.
In the final dimer dataset, the realized ranges are approximately $\delta \in [1.01\times 10^{-2},3.16\times 10^{-1}]$, $\gamma \in [1.00\times 10^{-3},9.95\times 10^{-2}]$, and $\kappa \in [1.00\times 10^{-2},9.99\times 10^{-1}]$.

Samples for which the nonlinear resonance solver fails to converge are discarded. We retain only outgoing damped resonances satisfying $\operatorname{Re}\omega>0$ and $\operatorname{Im}\omega<0$.

\subsection{Reference and reduced models}

For the single resonator, the asymptotic resonance is computed from the leading asymptotic eigenvalue of the volume integral operator. The reference resonance is computed from a more accurate numerical value of the same eigenvalue.

For the dimer, the reduced and reference eigenvalues are computed using two discretization levels of the two-resonator integral operator. The reduced model uses a coarser discretization, while the reference model uses a finer discretization. In the reported computations, the coarse dimer operator uses $25$ sample points per resonator and the reference operator uses $120$ sample points per resonator. The two resulting resonances are treated as the two dimer branches.

\subsection{Dimer discretization convergence}
\label{app:dimer_convergence}

For the dimer configuration, the reference resonances are obtained from a discretization of the two-resonator volume integral operator. Since the learned correction is trained against this numerical reference, it is important to check that the reference discretization is sufficiently stable. We therefore perform a convergence study with respect to the number $N$ of quadrature points per resonator.

For each value $N\in\{25,50,80,120,160,200\}$, we compute the two leading dimer resonances and compare them with the finest discretization $N_{\mathrm{ref}}=200$. The reported quantity is $|\omega_N-\omega_{N_{\mathrm{ref}}}|$, computed over the two branches and over a representative set of dimer configurations sampled from the same parameter ranges as in the main experiments.

\begin{table}[ht]
    \centering
    \begin{tabular}{lccc}
        \hline
        $N$ & Mean & Median & 90th percentile \\
        \hline
        $25$  & $3.56\times10^{-4}$ & $8.29\times10^{-5}$ & $1.18\times10^{-3}$ \\
        $50$  & $1.71\times10^{-4}$ & $3.88\times10^{-5}$ & $5.75\times10^{-4}$ \\
        $80$  & $9.68\times10^{-5}$ & $2.26\times10^{-5}$ & $3.21\times10^{-4}$ \\
        $120$ & $3.44\times10^{-5}$ & $6.97\times10^{-6}$ & $1.14\times10^{-4}$ \\
        $160$ & $1.69\times10^{-5}$ & $3.53\times10^{-6}$ & $5.53\times10^{-5}$ \\
        $200$ & $0$ & $0$ & $0$ \\
        \hline
    \end{tabular}
    \caption{Dimer discretization convergence. The errors are measured relative to the finest discretization $N_{\mathrm{ref}}=200$, using $|\omega_N-\omega_{N_{\mathrm{ref}}}|$, and are aggregated over both resonance branches.}
    \label{tab:dimer_discretization_convergence}
\end{table}

\begin{figure}[ht]
    \centering
    \includegraphics[width=0.5\textwidth]{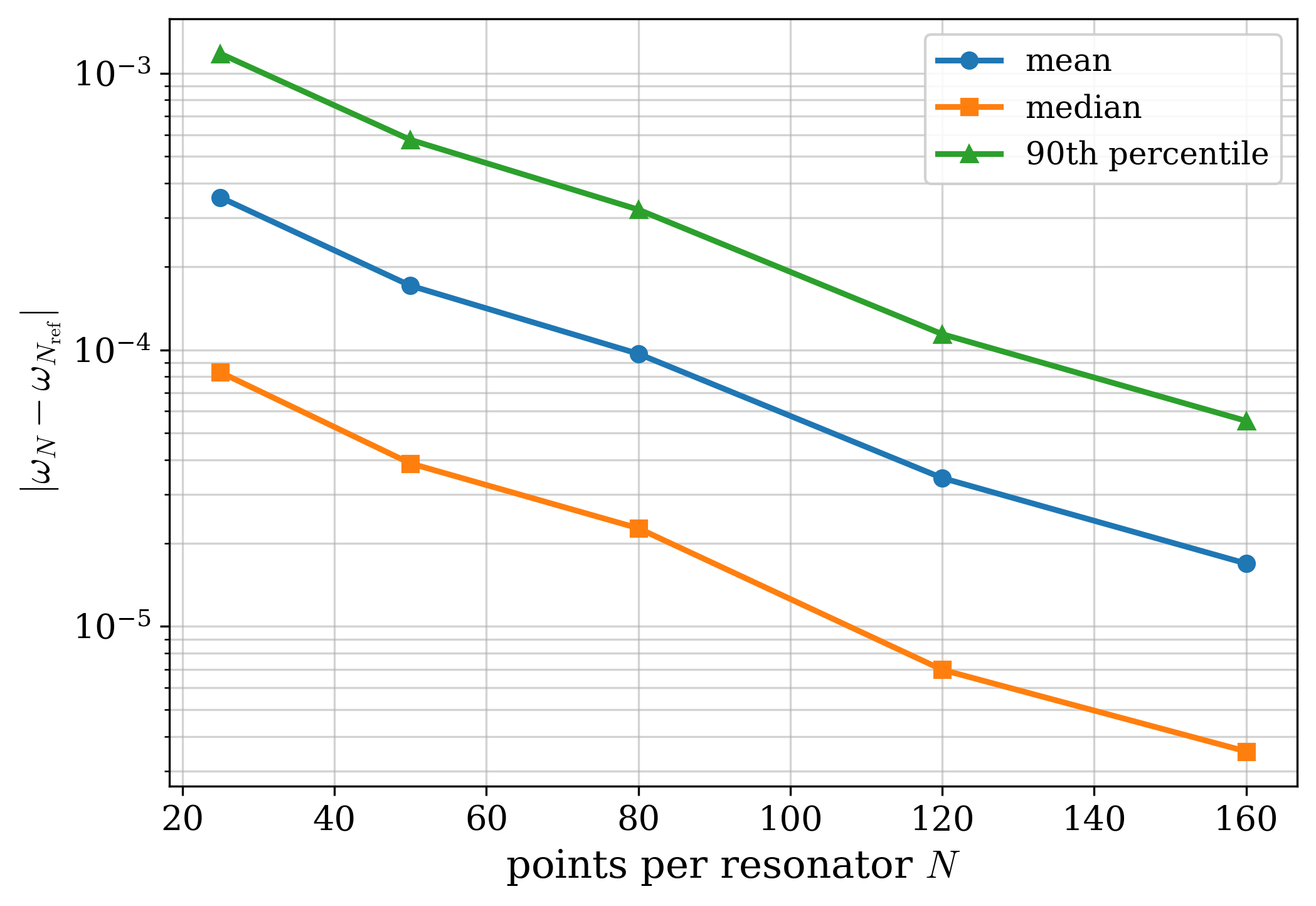}
    \caption{Dimer discretization convergence as a function of the number $N$ of quadrature points per resonator. The reference discretization is $N_{\mathrm{ref}}=200$.}
    \label{fig:dimer_discretization_convergence}
\end{figure}

The error decreases systematically as $N$ increases. In particular, the mean error drops from approximately $3.56\times10^{-4}$ at $N=25$ to approximately $1.69\times10^{-5}$ at $N=160$. This indicates that the reference discretization used in the learning experiments is stable enough for the observed learned corrections not to be merely an artefact of a coarse discretization.

\subsection{Train-test splits and validation protocols}

Unless otherwise stated, models are evaluated using an $80$--$20$ train-test split with fixed random seed. This random split is used for the in-distribution tests. For out-of-distribution validation, the split is performed by holding out a challenging region of parameter space. For the single resonator, we test extrapolation to large $\delta$ and large $\gamma$. For the dimer, we test extrapolation to large $\delta$ and small $\kappa$. The holdout threshold is taken at the corresponding $80\%$ quantile, except for the small-$\kappa$ dimer test, where the lower $20\%$ of $\kappa$-values is held out. For a test sample $j$, the relative error $E_j$ and the reported mean relative error $\overline{E}$ are given in \eqref{eq:to be said}.
For the dimer, the average is taken over both resonance branches.

For the main prediction figures, the single-resonator dataset contains $10\,000$ samples, split into $8\,000$ training samples and $2\,000$ test samples. The dimer dataset contains $3\,000$ configurations, split into $2\,400$ training configurations and $600$ test configurations. For the large-data noisy-training experiments, the single-resonator dataset contains $12\,500$ samples, split into $10\,000$ training samples and $2\,500$ test samples. The dimer noisy-training experiment uses $5\,000$ training configurations and $1\,250$ test configurations.

\subsection{Ridge regression}

The Ridge model is trained on the residual
\begin{align*}
    \Delta\omega =
    \omega_{\mathrm{ref}}-\omega_{\mathrm{asymp}},
\end{align*}
rather than on the full resonance. The real and imaginary parts of the residual are learned simultaneously.

The Ridge objective is
\begin{align*}
    \min_{W,b}
    \|Y-XW-b\|_2^2
    +
    \alpha_R \|W\|_2^2,
\end{align*}
with $\alpha_R = 10^{-8}$.
Before training, all input features are standardized using the empirical mean and standard deviation of the training set. The same transformation is then applied to the test set. We do not tune $\alpha_R$ by cross-validation; the objective is to test whether a simple and stable linear residual model can exploit the asymptotic feature structure, rather than to optimize a large hyperparameter search.
For the calibrated next-order proxy, we use a smaller Ridge parameter, $\alpha_R = 10^{-12}$,
since this step is used as a least-squares calibration of the candidate asymptotic feature family.

\subsection{Feature sets}

The raw feature sets contain logarithmic transforms of the sampled parameters, such as $\log\delta$, $\log\gamma$ and, in the dimer case, $\log\kappa$.
The physics-informed feature sets include combinations motivated by the asymptotic scaling, for example $\delta^2$, $\delta^2\log\delta$, $\gamma\delta^2$.

The frequency-dependent next-order terms are motivated by the scale \eqref{eq:scale}. In the implementation, the independent logarithmic features are therefore taken to be $\delta^2 q_{\mathrm{asymp}}^2\log\delta$, $\delta^2 q_{\mathrm{asymp}}^2\log q_{\mathrm{asymp}}$, together with $\delta^2 q_{\mathrm{asymp}}^2$.

In the dimer case, the frequency-dependent features are evaluated branchwise using the corresponding asymptotic resonance. Additional coupling features include the asymptotic eigenvalues and the asymptotic branch splitting.

\subsection{Symbolic regression}

Symbolic regression is used to obtain compact closed-form corrections for the residual. Separate formulas are fitted for the real and imaginary parts. In the dimer case, symbolic regression is applied separately to each branch.

For the clean-data PySR runs used to obtain the symbolic formulas reported in the main text, the allowed binary operators are $+,\quad -,\quad \times$,
and no unary operators are used. The model selection criterion is PySR's \texttt{best} rule. The maximum expression size is set to $20$, with $20$ populations and $500$ iterations. Fixed random seeds are used.

For noisy symbolic rediscovery tests, PySR is rerun with the operator set $+,\quad -,\quad \times,\quad /$,
maximum expression size $20$, population size $40$, deterministic serial execution, and $100$ iterations unless otherwise specified. These runs are used only for robustness diagnostics. The compact symbolic formulas reported in the main text are the clean-data formulas.

\subsection{Assisted symbolic regression}

In the assisted symbolic-regression experiments, a calibrated next-order proxy is first fitted using the frequency-dependent feature family described above. We then compare two procedures:

\begin{enumerate}
    \item symbolic regression is applied to the remaining residual after subtracting the proxy;
    \item the proxy is supplied as an additional input feature to symbolic regression.
\end{enumerate}

These experiments test whether explicit next-order asymptotic information improves the symbolic search.

\subsection{Noise experiments}

For the noisy-training experiments, the training resonance is perturbed according to
\begin{align*}
    \omega_{\mathrm{ref}}^{\mathrm{noisy}}
    =
    \omega_{\mathrm{ref}}
    +
    \sigma |\omega_{\mathrm{ref}}|
    \frac{\eta_1+i\eta_2}{\sqrt{2}},
\end{align*}
where $\eta_1$ and $\eta_2$ are independent standard Gaussian random variables. The models are trained on
\begin{align*}
    \omega_{\mathrm{ref}}^{\mathrm{noisy}}
    -
    \omega_{\mathrm{asymp}},
\end{align*}
and evaluated on a clean held-out test set.

The noise levels are $\sigma \in \{0,10^{-4},10^{-3},10^{-2}\}$.
For the Ridge and proxy-based noise experiments, each noise level is averaged over $20$ independent noise realizations. For the symbolic rediscovery experiments, PySR is rerun independently for each noise level.

In the clean-versus-noise-trained symbolic experiment, the test resonances are also perturbed. This experiment compares a symbolic formula discovered from clean data with a symbolic formula rediscovered from noisy training data, both evaluated against noisy test resonances.

\bibliographystyle{abbrv}
\bibliography{references}{}

\end{document}